\documentclass[11pt]{amsart}
\usepackage{times}      
\usepackage{latexsym}   
\usepackage{amssymb}    
\usepackage{amsmath}    
\usepackage{amsthm}     
\usepackage{mathrsfs}   
\usepackage{euscript}   
\usepackage{epsfig,subfigure,citesort}             
\usepackage[all]{xy}    
\newcommand{\cross}{\times}

\newcommand{\ra}{\rightarrow}
\newcommand{\R}{\mathbb{R}}

\newcommand{\Hy}{\mathbb{H}^1}

\renewcommand{\epsilon}{\varepsilon}
\newcommand{\nin}{\not\in}
\renewcommand{\qed}{\square}
\renewcommand{\"}[1]{\mathaccent 127 #1}
\renewcommand{\phi}{\varphi}

\newcommand{\cc}{d_{cc}}

\renewcommand{\"}[1]{\mathaccent 127 #1}
\renewcommand{\phi}{\varphi}

\newcommand{\Ha}{\mathscr{H}}
\newcommand{\ip}{<\cdot,\cdot>}
\newcommand{\pf}{\noindent {\em Proof: }}
 
\newtheorem{Pro}{Proposition}[section]
\newtheorem{Lem}[Pro]{Lemma}

\newtheorem{Thm}[Pro]{Theorem}
\newtheorem{MThm}{Theorem}

\newtheorem{Rem}{Remark}
\newtheorem{Def}[Pro]{Definition}

\newtheorem{Cor}[Pro]{Corollary}

\newtheorem{eg}{Example}

\newcommand{\Om}{\Omega}

\newcommand{\pb}{\overline{p}}
\newcommand{\qb}{\overline{q}}
\newcommand{\pte}{\tilde{p}_\epsilon}
\newcommand{\qte}{\tilde{q}_\epsilon}

\newcommand{\nux}{{\boldsymbol{\nu}_X}}
\newcommand{\nuxp}{\nux^\perp}

\begin{document}

\title{H-minimal graphs of low regularity in $\Hy$}
\author{Scott D. Pauls}
\address{Dartmouth College, Hanover, NH 03755}
\thanks{The author is partially supported by NSF grant DMS-0306752}
\email{scott.pauls@dartmouth.edu}

\begin{abstract}
In this paper we investigate H-minimal graphs of lower regularity.  We
show that noncharactersitic $C^1$ H-minimal graphs whose components of the unit horizontal Gauss map are in $W^{1,1}$ are ruled surfaces with $C^2$ seed
curves.  In a different direction, we investigate ways in which patches of $C^1$ H-minimal graphs
can be glued together to form continuous piecewise $C^1$ H-minimal
surfaces.

We apply these description of H-minimal graphs to the question of the
existence of smooth solutions to the Dirichlet problem with smooth
data.  We find a necessary condition for the existence
of smooth solutions and produce examples where the conditions are
satisfied and where they fail.  In particular we illustrate the failure
of the smoothness of the data to force smoothness of the solution to
the Dirichlet problem by producing a class of
smooth curves whoses H-minimal spanning graphs cannot be $C^2$.

\end{abstract}

\maketitle
\section{Introduction}
In this paper, we further investigate the properties of H-minimal
surfaces in the Heisenberg group with a focus on the regularity of H-minimal surfaces that
satisfy Dirichlet boundary conditions. 

The study of H-minimal surface
was introduced in the foundational paper of Garofalo and Nhieu
(\cite{GarNh}) where they showed the existence of H-minimal surfaces
of bounded variation that satisfied certain boundary conditions.
Expanding on these results, several authors extended the investigation
showing different properties and constructions of H-minimal surfaces
in various settings
(see, for example, \cite{Pauls:minimal,DGN,GP,CHMY,CH,RR,Cole,HP2,HP,BC}).  Recently, N. Garofalo and the author (\cite{GP}) gave a characterization of
$C^2$ H-minimal surfaces used to investigate an analogue of the Bernstein
problem in the Heisenberg group.  A different approach to the study of
analogues of the Bernstein problem 
was completed by Cheng, Hwang, Malchiodi and Yang in \cite{CHMY} and \cite{CH}.  We note that there is some overlap between the results in \cite{GP} and those of
\cite{CHMY} and \cite{CH} but that the techniques are independent.  In
particular, both \cite{GP} and \cite{CHMY} make the observation that
$C^2$ H-minimal surfaces are ruled surfaces but analyze them using
different tools (in fact, \cite{CHMY} uses the machinery of
pseudohermitian geometry and hence many of their results apply to a larger
class of Carnot-Carath\'eodory spaces).  Using the machinery of
\cite{CHMY}, two of the authors 
classify properly embedded H-minimal surfaces in the Heisenberg
group in \cite{CH} while \cite{GP} gives a geometric description of the properties
of embedded H-minimal surfaces which are graphs over some plane.  Again, the results overlap in some
respects, but the techniques are independent.

With respect to the discussion in this paper, we will use the tools
developed in \cite{GP}.  For the purposes of this paper, the two most important
theorems from \cite{GP} are: 

\begin{Thm}\label{mrep} Let $k \ge 2$.  
A noncharacteristic patch of a $C^k$ surface $S \subset \Hy$ of the type
\[
S\ =\ \{(x,y,t)\in \mathbb H^1 \mid (x,y)\in \Om\ ,\ t = h(x,y)\}\ ,
\]
where $h : \Om \to \R$ is a $C^k$ function over an open domain $\Om$ in the
$xy$-plane, is an $H$-minimal surface if and only if for every $p \in
S$, there exists a neighborhood $U$ of $p$ so that $U$ can be
parameterized by 
\begin{equation}
(s,r)\ \to\ (\gamma_1(s)+r\gamma_2'(s),\gamma_2(s)-r\gamma_1'(s),
h(s,r))\ ,
\end{equation}
where
\begin{equation}
h(s,r)\ =\ h_0(s)\ -\ \frac{r}{2} <\gamma(s), \gamma'(s)>\ .
\end{equation}
and \[ \gamma \in C^{k+1}, \; h_0 \in C^k\]
Thus, to specify such a patch of smooth $H$-minimal surface, one must specify a
single curve in $\Hy$ determined by a curve in the plane,
$\gamma(s)$, parameterized by arc-length, and an initial height function $h_0(s)$.
\end{Thm}

The curve $\gamma(s)$ in the theorem is called a {\bf seed curve} and
determines almost all of the behavior of neighborhood $U$.  Indeed, under the assumption of at least $C^2$ smoothness, we have:

\begin{Thm}\label{cleanup}Let $S\subset \mathbb H^1$ be a $C^2$
  connected, open, complete and embedded $H$-minimal surface.
  Then, either $S$ is a vertical plane, or $S$ is determined by a generalized 
  seed curve $\Gamma = \{(\gamma_1^i(s),\gamma_2^i(s),h_0^i(s))\}$.  
\end{Thm}

A generalized seed curve is a collection of seed curves, height
functions and patching data which, taken togeher, give a description
of a single curve in $\mathbb{H}^1$.  In other words, for such
H-minimal surfaces, a single curve determines the entire surface. 

As all H-minimal surfaces have locally finite perimeter (i.e. they are
X-Caccioppoli sets), we turn now to the work of Franchi, Serra Cassano and Serapioni
(\cite{FSSC2}) and  recall the following theorem:

\begin{Thm}\label{decomp}
Let $E \subset \Hy$ be a $X$-Caccioppoli set, then the reduced boundary of $E$, $\partial^*_X E$, is $X$-rectifiable, i.e.,  
\[
\partial^*_X E\ =\ N\ \cup\ \bigcup_{j=1}^\infty K_j\ ,
\]
where $\mathcal H_{cc}^{3}(N) = 0$, and $K_j$ is a compact subset of a non-characteristic hypersurface $S_j$ of class $C^1_\mathbb{H}$. Moreover, one has for any $g\in K_j$ and every $\xi \in T_{H,g} S_j$
\[
<\nu^E_X(g),\xi>\ =\ 0\ ,
\]
where $\nu_E(g)$ denotes the generalized horizontal outer normal to
$E$ in $g$, $T_{H,g} S_j$ indicates the non-characteristic plane
orthogonal to the horizontal normal to $S_j$ in $g$, $\mathcal
H_{cc}^{3}$ is the $3$-dimensional Hausdorff measure in $\Hy$
constructed with respect to the CC distance and $C^1_\mathbb{H}$ is the
space of functions which are horizontally continuously differentiable,
i.e. $X_1f,X_2f$ exist and are continuous.
\end{Thm}

The reduced boundary, $\partial^*_X$, is the set of boundary points
where the unit horizontal Gauss map is well-defined (see the next
section for a precise definition).  For the discussions of this paper,
it is important to note that the reduced boundary is a full measure
subset of the boundary.  
The main point of this theorem is that H-minimal surfaces can be decomposed into a set of
$\Ha^3_{cc}$-measure zero and a union of $C^1_\mathbb{H}$ sets.  As we will restrict
ourselves to investigating graphs over the xy-plane, we remind the
reader that a $C^1_\mathbb{H}$ graph is $C^1$.  This, of course,
leaves a gap - the pieces given by theorem \ref{decomp} are $C^1$
while theorem \ref{cleanup} applies only to $C^2$ surfaces.  The first
goal of this paper is to partially bridge the gap between the two theorems.
\begin{MThm}\label{c1surf}  If $S$ is an open $C^1$ H-minimal graph over a domain
  $\Omega \subset \R^2$ with no
  characteristic points and unit horizontal Gauss map $\nux$ whose components are in $W^{1,1}(\Omega)$, then
  the integral curves of $\nux^\perp$ are straight lines and $S$ can
  be locally parameterized by
\begin{equation}\label{Mpara}
(s,r)\ \to\ (\gamma_1(s)+r\gamma_2'(s),\gamma_2(s)-r\gamma_1'(s),
h(s,r))\ ,
\end{equation}
where
\begin{equation}\label{Mpara2}
h(s,r)\ =\ h_0(s)\ -\ \frac{r}{2} <\gamma(s), \gamma'(s)>\ 
\end{equation}
and $\gamma$ is an integral curve of $\nux$.  Moreover, if there
exists $\epsilon >0$ so that $d(\Omega,\gamma)
>\epsilon >0$ then $h_0(s) \in C^1$ and $\gamma(s) \in C^2$.
\end{MThm}
In this theorem, $d(\Omega,\gamma)$ is a measure of the ``horizontal
thickness'' of the set $\Omega$ (see definition \ref{hthick} for a precise statement).

This theorem is shown in a series of steps.  First, we show that the
weak directional derivative of $\nux$ in the direction of $\nuxp$ is
zero.  This is enough to show that the integral curves of $\nuxp$ are
lines.  Second, forming $\gamma(s)$ as the integral curve of $\nux$, a
geometric argument shows that $\gamma'(s)$ is Lipschitz.  Coupled with
a further estimate, this shows that $\gamma''(s)$ exits and is continuous.
Applying arguments similar to those in \cite{GP} yields the
representation given in the theorem.  

Combining this theorem with theorem \ref{decomp} of Franchi, Serapioni and Serra Cassano yields:

\begin{MThm}\label{sum} If $S$ is an H-minimal graph then 
\[S = N \cup \bigcup_{i=1}^\infty K_i\]
Where $N$ is a set of $\Ha^3_{cc}$-measure zero and each $K_i$ is a compact piece of a $C^1$ H-minimal graph.  For each $i$, if $K_i$ is open and its unit horizontal Gauss map has components in $W^{1,1}$ then $K_i$ may be locally parameterized by equations \eqref{Mpara} and
\eqref{Mpara2} with $\gamma \in C^2$ and $h_0 \in C^1$. 
\end{MThm}

While this theorem recovers the characterization of such H-minimal
surfaces as unions of ruled surfaces in the Heisenberg group, it still leaves a gap between the baseline results of Franchi, Serra
Cassano and Serapioni and theorem
\ref{mrep}.  In particular, theorem \ref{decomp} allows that the $C^1$
pieces may be glued together in nonsmooth ways.  We find that this can
happen:

\begin{MThm}\label{pwc1}
Suppose $S_1$ and $S_2$ are subsets of $C^1$ H-minimal graphs
  with no characteristic points, each parameterized by a single seed
  curve and height function, defined over closed sets $\Omega_1,\Omega_2
  \subset \R^2$ with open interior, $C= \Omega_1
  \cap \Omega_2$ a $C^1$ curve, $\partial \Omega_i \in C^1$ and $d(\Omega_i)>0$ for $i=1,2$.  Moreover, let $\nu_1=(\pb_1,\qb_1)$
  and $\nu_2=(\pb_2,\qb_2)$ be the respective unit horizontal Gauss
  maps.  Then, $S_1 \cup S_2$ is an H-minimal graph if and only if
  $(\nu_1-\nu_2)|_C$ is tangent to $C$ almost everywhere. 
\end{MThm}
  This provides one way in which a
continuous, piecewise $C^1$ H-minimal graph is constructed.

This brings forward an obvious question:

\vspace{.2in}

\noindent
{\em In standard minimal surface theory, the solutions to the minimal
  surface equation subject to Dirichlet boundary conditions gain
  additional regularity from the regularity of the boundary.  Do H-minimal surfaces have a similar property?}

\vspace{.2in}

We devote the remaining part of the paper to exploring this question.  First, we
examine some of the best behaved H-minimal surfaces, those that are
minimal in Riemannian approximators of $\Hy$ as well as H-minimal.  We
call these {\bf persistent minimal surfaces} and classify them:
  
\begin{MThm}  The persistent H-minimal graphs fall into two
  categories:
\begin{enumerate}
\item $S$ is given by $(x,y,u(x,y))$ where 
\begin{equation*}
\begin{split}
u(x,y)&=\frac{m}{1+m^2}(x-x_0)^2+\frac{m^2-1}{m^2+1}(x-x_0)(y-y_0) -
    \frac{m}{1+m^2}(y-y_0)^2+\\
   & \frac{a}{\sqrt{1+m^2}}(x-x_0)+\frac{am}{\sqrt{1+m^2}}(y-y_0)+b
\end{split}
\end{equation*}
for $m,a,b,x_0,y_0 \in \R$.
\item $S$, given in cylindrical coordinates is \[(\rho
  \cos(\theta),\rho \sin (\theta), a \theta + b)\] for $a,b \in \R$
\end{enumerate}
\end{MThm} 

These surfaces give examples of the best possible case - they are
$C^\infty$ spanning surfaces.  Second, in section 7, we consider the
question of the existence of smooth minimal spanning surfaces.  For a
fixed smooth closed curve, we focus on finding the $C^1$
smooth {\bf ruled H-minimal spanning graphs}, those $C^1$ H-minimal graphs which are ruled surfaces away from the characteristic locus which satisfy the additional condition that the rules may be extended over the characteristic locus (as straight lines).  We note that by work in \cite{GP}, all $C^2$ H-minimal graphs satisfy this condition.  However, using techniques similar to those of theorem
\ref{pwc1}, one can construct $C^1$ minimal graphs that are the union
of two ruled surfaces along a common characteristic locus and do not satisfy this condition (see section 4, example 4).  For
simplicity, we will ignore this type of construction.  Taking the characterization of
H-minimal surfaces as ruled surfaces, we create a necessary condition for a smooth closed curve which is the graph over
a curve in the xy-plane to be spanned by a $C^1$ ruled H-minimal graph:

\vspace{.25in}
\noindent
{\bf Existence Criteria:}  {\em Given a closed smooth curve $c(\theta)=(c_1(\theta),c_2(\theta),c_3(\theta))$ which is a
graph over a curve in the xy-plane, if $c$ is spanned by a $C^1$
ruled H-minimal graph then there exists a monotone $C^1$
function $\phi :S^1 \ra S^1$ with $\phi(\theta) \in A(\theta)$.  }

\vspace{.25in}

In this statement, $A(\theta)$ is the set of points on $c$ that are
accessible from $c(\theta)$ via a rule of an H-minimal surface:

\[A(\theta)=\left \{\theta_0 |
  c_3(\theta_0)-c_3(\theta)-\frac{1}{2}c_1(\theta_0)c_2(\theta)
  +\frac{1}{2}c_1(\theta)c_2(\theta_0) = 0 \right \}\]

The examples in this section show curves that satisfy the criteria and
curves that exhibit an obstruction.  We also discuss the genericity of
these classes.  Finally, we show that there are many curves, $c$,
which do not have smooth ruled H-minimal spanning graphs.  This provides an
upper bound on the regularity of the solution to the Plateau Problem
for these curves:  the solution to the Plateau Problem cannot be $C^2$.

\begin{MThm}
 Suppose $c$ is a $C^1$ curve with no Legendrian points
  which is spanned by a $C^1$ smooth ruled H-minimal graph, $S$.  Then there
  exists an interval, $I$, so that $c(I)$ is contained in a plane. 
\end{MThm}

\begin{Cor}
If $c$ is a smooth curve with no Legendrian points and no
  portion of 
  $c$ is contained in a plane then an H-minimal surface spanning
  $c$ cannot be a $C^1$ ruled H-minimal surface.  
\end{Cor} 

These different examples show that solutions to the Dirichlet problem
and the Plateau Problem may not have any specified regularity.  In
particular, the persistent H-minimal graphs show that some curves
have $C^\infty$ solution to the Plateau Problem while the subsequent
examples show instances where $C^\infty$ curves may not have solution
to the Dirichlet problem of high regularity.  Indeed, the last set of
examples show that for certain totally non-Legendrian curves, the
graphical solutions to the Dirichlet (and hence the Plateau) problem are
neccesarily at most $C^1$ but cannot be ruled surfaces.  A
consequences of this is that these surfaces must have unresolvable
discontinuities in their unit horizontal Gauss maps.  

The author would like to thank J. H> Cheng, J. F. Hwanf and P. Yang for pointing out an error in an earlier version of this paper.  The author would also like to thank the referee for many helpful comments and suggestions.  

\section{Definitions and Notation}\label{defs}
Throughout this paper, we restrict our attention to the topologically
three dimensional Heisenberg group, $\Hy$.  For convienience, we
represent $\Hy$ via an identification with $\R^3$.  Considering $\R^3$
with its usual coordinates labeled as $\{x,y,t\}$, we define the
following vector fields:

\begin{equation*}
\begin{split}
X_1 &= \frac{\partial}{\partial x} - \frac{y}{2}
\frac{\partial}{\partial t}\\
X_2 &= \frac{\partial}{\partial y} + \frac{x}{2}
\frac{\partial}{\partial t}\\
T &= \frac{\partial}{\partial t}
\end{split}
\end{equation*}

The vector fields $\{X_1,X_2,T\}$ form a basis for the Lie algebra of
$\Hy$ at any point $(x,y,t)$.  
Note that, via the exponential map at the origin, we
identify $\Hy$ with $\R^3$ using these coordinates,  denoting the point $e^{\alpha
  X_1+\beta X_2+\gamma T}$ by $(\alpha,\beta,\gamma)$.  For the purposes
of this paper, we define a left invariant inner product on $\Hy$, $\ip$,
which makes $\{X_1,X_2,T\}$ an orthonormal basis at each point.
Notice that at each point, $[X_1,X_2] =
T$ and hence $\{X_1, X_2\}$ is a bracket generating set for $\Hy$.  
We define a subbundle on $\Hy$, called the horizontal subbundle of $\Hy$,
by 
\[H\Hy = \{(x,y,t,w) \in \Hy \cross \R^3 | w \in span\{X_1,X_2\}\}\]
The single
nontrivial bracket relation yields the following multiplication law
via the Campbell-Baker-Hausdorff formula:

\[(a,b,c)(\alpha,\beta,\gamma) = \left (a+\alpha,b+\beta,c+\gamma +
\frac{1}{2}(a\beta - \alpha b)\right )\]

To define the Carnot-Carath\'eodory metric on $\Hy$, we construct a
path metric.  Letting $\mathscr{A}(m,n)$ be the set of all absolutely
continuous paths in $\Hy$ so that $\gamma(0)=m, \gamma(1)=n$ and
$\gamma'(t) \in  H_{\gamma(t)}\Hy \text{ when $\gamma'(t)$ exists}$ we define:

\[ \cc(m,n) = \inf_{\gamma \in \mathscr{A}(m,n)}  \left \{\int_I <\gamma'(t),\gamma'(t)>^\frac{1}{2} \right \} \]

Note that, since $\ip$ is left invariant, so is $\cc$.  Moreover,
$\cc$ admits a homothety at each point $(x,y,t)$:

\[ h_s(x,y,t) = (sx,sy,s^2t)\]

whereby 
\[\cc(h_sm, h_sn)=s\cc(m,n)\]

We denote by $\Ha^k_{cc}$ the k-dimensional spherical Hausdorff
measure constructed from $\cc$.

\begin{Def}
The {\bf horizontal gradient operator} is 
\[\nabla_0 = (X_1,X_2)\]
Hence, 
\[\nabla_0 f = (X_1 f) \; X_1 + (X_2 f) \; X_2\]
The {\bf horizontal divergence} of a vector field $V=v_1 \;X_1+v_2\;X_2$ is
\[div_0 V=\nabla_0 \cdot V = X_1 v_1+X_2 v_2\]
\end{Def}

In addition to the three dimensional Hausdorff measure, we recall
the perimeter measure introduced independently by Capogna, Danielli and Garofalo
(\cite{CDG,CDG2}) and Franchi, Gallot and Wheeden (\cite{FGW}).

\begin{Def}  Let $\Omega$ be an open subset of $\Hy$.  We say that $f:
  \Omega \ra \R$ is of bounded variation (i.e. $f \in BV_\Hy(\Omega)$)
  if 
\begin{itemize}
\item $f \in L^1(\Omega)$
\item \[\sup \left \{ \int f \; div_0 V \; dh | V \in C^1_0(\Omega,H\Hy),
  |V| \le 1\right \} < \infty\]
\end{itemize}
\end{Def}
We define $BV_{\Hy,loc}(\Omega)$ analogously.

\begin{Def} We say that $E \subset \Hy$ is an {\bf X-Caccioppoli set} if the
  characteristic function of $E$, $\chi_E$ if it is in
  $BV_{\Hy,loc}$. The measure $|\nabla_0 \chi_E|$ is called the
  perimeter measure and will be denoted by $\mathscr{P}$.
\end{Def}
We recall that (see \cite{GarNh,FSSC2}) if $\partial E$ is a
smooth surface given by $t=u(x,y)$, then $\mathscr{P}$ is
mutually absolutely continuous with $\Ha^3_{cc}$.  Moreover, up to a
choice of constant, $\Ha^3_{cc}(\partial E)$ is given by
\[ \int_{\partial E} \sqrt{(X_1(t-u(x,y)))^2+(X_2 (t-u(x,y)))^2} \; dA \]

Again from \cite{FSSC2}, we recall the definition of the generalized
horizontal normal:
\begin{Def}
There exists a $\mathscr{P}$ measurable section $\nu_E$ of $H\Hy$ such
that
\[ - \int_E div_0 \phi\; dh = \int_\Hy <\nu_E, \phi> d\mathscr{P}\]
for all $\phi \in C^\infty_0(\Omega,H\Hy)$, $|\nu_E(p)|=1$ for
$\mathscr{P}$ a.e. $p \in \Hy$.
\end{Def}

We next recall the definition of the reduced boundary:
\begin{Def} Let $E$ be an X-Caccioppoli set.  Let $U(p,r)$ be the
  open ball of radius $r$ and center $p$.  A point $p$ is in the
  {\bf reduced boundary of E}, $p \in \partial^*_{\Hy} E$, if 
\begin{enumerate}
\item $\mathscr{P}(U(p,r))>0$ for any $r>0$
\item \[\lim_{r \ra 0} \frac{1}{\mathscr{P}(U(p,r))}\int_{U(p,r)} \nu_E \;d \mathscr{P}\] exists
\item \[\left |\lim_{r \ra 0} \frac{1}{\mathscr{P}(U(p,r))}\int_{U(p,r)}
  \nu_E \;d \mathscr{P}\right| =1\]
\end{enumerate}
\end{Def}
We note that lemma 7.3 in \cite{FSSC2} ensures that $\partial^*_{\Hy}E$
has full $\mathscr{P}$ measure in $\partial E$.  

In this paper, we will be examining smooth graphs over the xy-plane in $\Hy$
by which we mean surfaces which can be represented as $t=u(x,y)$ using
the coordinates described above.  As shown in \cite{Pauls:minimal} and
\cite{DGN}, H-minimal surfaces are critical points of an area
functional based on the {\bf horizontal Gauss map} of the suface
$t=u(x,y)$.  The horizontal Gauss map is the projection of
the Riemannian normal of the surface to the horizontal bundle:
\begin{equation*}
\begin{split}
 G:S &\ra HS\\
G(x,y,u(x,y))&= (X_1 (t-u(x,y)))\;X_1+ (X_2 (t-u(x,y)))\;X_2
\end{split}
\end{equation*}

We give classically inspired names to these horizontal derivatives of
$f$, letting 
\[p= X_1 (t-u(x,y))\]
\[q= X_2 (t-u(x,y))\]

In this paper, the {\bf unit horizontal Gauss map} plays a crucial role and
so we define the unit horizontal Gauss map by 
\[\nu_X = \pb X_1 + \qb X_2\]
where $\pb = \frac{p}{\sqrt{p^2+q^2}}$ and  $\qb=
\frac{q}{\sqrt{p^2+q^2}}$.  Notice that $\nux$ has a limited domain
and is not defined at points where both $p$ and $q$ are zero.  Such
points are called {\bf characteristic points} and play an important role in
the study of surfaces in Carnot-Carath\'eodory spaces.  

In this paper, we consider surfaces which are graphs over the
xy-plane.  In other words, the set $E$ in the previous theorem is
given as 
\[\{(x,y,t)| t < u(x,y)\}\]
Thus, the hypersurface $\partial E$ would be given as $t-u(x,y)=0$.
The function $t-u(x,y)$ is horizontally continuously
differentiable if and only if $u$ is continuously differentiable.  

With this notation in place, we next review the characterization of
smooth noncharacteristic area minimizing graphs by an appropriate partial differential
equation via the first variation of the energy.  The first variation
formula has been explored in a variety of settings by a number of authors including the
aforementioned paper of Cheng, Hwang, Malchiodi and Yang
(\cite{CHMY}), Danielli, Garofalo and Nhieu (\cite{DGN}), the author (\cite{Pauls:minimal}), Bonk and Capogna
(\cite{BC}), Ritor\'e and Rosales (\cite{RR}).  For the convienience of the reader, we recall the 
derivation of the equation here.  First, the energy
integral we use for the variational setup is: 

\[E(u) = \int_\Omega \sqrt{p^2+q^2} \]

where $t=u(x,y)$ defines the graph in question over a domain
$\Omega$ and has at least two weak derivatives.

Second, we consider a variation in the $t$ direction by a function
$\phi(x,y) \in C_0^\infty(\Omega)$.  Then, 

\[E(\epsilon)=E(u+\epsilon \phi)= \int_\Omega   \sqrt{\left(u_x+\frac{y}{2}+\epsilon
  \phi_x\right)^2+\left(u_y-\frac{x}{2}+\epsilon \phi_y\right)^2}\]

Abusing notation, we let $p=u_x+\frac{y}{2}$ and $q=u_y-\frac{x}{2}$.
Thus, differentiating with respect to $\epsilon$ twice and evaluating at zero, we have:

\[ E'(0)=\int_\Omega \frac{p\phi_x+q\phi_y}{\sqrt{p^2+q^2}}\]

and 

\[ E''(0) = \int_\Omega \frac{(q\phi_x-p\phi_y)^2}{(p^2+q^2)^\frac{3}{2}} =
\int_\Omega \frac{(\nabla \phi \cdot G^\perp)^2}{|G|^3}\]

In the last equation, we use the convention that if $v$ is the vector
given by coordinates $(a,b)$ then $v^\perp$ is given by $(b,-a)$.
This convention will be used throughout the paper. Note that the
integrand of the second integral is nonnegative and is strictly
positive if $\nabla \phi$ is not parallel to the vector $G$.  Thus, to check
if a given solution to the Euler-Lagrange equation is a local minimum
of area (with respect to this type of perturbation), one must only check it against variations in this direction.

\begin{Lem} Let $u:\Omega \ra \R$, $u \in C^2$, be a critical point of the energy
  functional.  Then, for all test functions $\phi \in
  C^\infty_0(\Omega)$, $E''(0)>0$.
\end{Lem}

\pf  Notice first that if $\nabla \phi \cdot G^\perp$ is not identically
zero on a set of full measure, then since the integrand is always positive, the result
follows.  If $\nabla \phi$ points in the same direction as $G$, we now
verify that for such a perturbation, $E'' > 0$.  In this case, let  $\beta$ be a
function so that $\nabla \phi =
\beta(x,y) G=(\beta G_1,\beta G_2)$.  Then, $\beta G$ is the gradient
of the $C^\infty_0(\Omega)$ function $\phi$ so,
\[\phi_{xy}-\phi_{yx} = \beta_yG_1+\beta G_{1,y} -\beta_xG_2-\beta G_{2,x}
=0\]

But, 
\[\beta G_{1,y} -\beta G_{2,x}= \beta(u_{xy}-u_{yx}-1)=-\beta\]

So, we have
\begin{equation}\label{aa}
 -\nabla \beta \cdot G^\perp-\beta =0
\end{equation}
By theorem \ref{mrep} in the introducton, we see that the integral curves of $G^\perp$ are straight lines (for a more detailed dicussion of this fact, see \cite{GP}).  Thus, \eqref{aa} implies that when $\beta$ is restricted to such a straight line, we have $\beta'=-\beta$ and hence, $\beta = C
e^{-t}$ where $t$ is the parameter along the integral curve.  However,
as $\phi$ is compactly supported in
$\Omega$, $\beta$ must tend to zero towards the boundary of $\Omega$.
This is a contradiction of the existence of a
compactly supported normal variation $\phi$ with gradient pointing in the
same direction as $G$.  $\qed$.

Thus local minima of this area functional appear as solutions of
the following partial differential equation:

\begin{equation}\label{MSE}
 X_1 \pb + X_2 \qb =0
\end{equation}

This equation says simply that the unit horizontal Gauss
map is (horizontally) divergence free:

\[ div_0 \; \nux =0\] 
In this paper, we will also allow solutions that are only weak
solutions to this equation.  In
section \ref{glue}, we discuss a condition under which a piecewise $C^1$ graph
can satisfy this equation weakly, but not strongly.  

\begin{Rem}
In some of the references given above, there are versions of these first and second variation
equations for more generally defined surfaces. For example,in \cite{DGN},
the authors give these formulae for surfaces defined implicitly asa
$\phi(x,y,t)=0$.  We also point out that without our restriction to
graphs (and perterbations that remain graphs) the second variation
formula does not necessarily ensure that the critical points are local
minima.  Again, see the examples in \cite{DGN}.
\end{Rem}

In \cite{DGN}, Danielli, Garofalo and Nhieu introduce the notation of
H-mean curvature, which is used to define H-minimal surfaces in both
\cite{DGN} and \cite{GP}.  We recall slight variations of these
definitions here using the notation above.  
\begin{Def}The {\bf H-mean curvature} of $S$ at noncharacterstic
  points of $S$ is defined by 
\[ H =  X_1 \pb + X_2 \qb \]
If $x_0 \in \Sigma$, then 
\[H(x_0)= \lim_{x \ra x_0, x \nin \Sigma} H(x)\]
provided that the limit exists, finite or infinite.  If the limit
does not exist, the H-mean curvature is not defined at such points.
\end{Def}

This definition differs from that in \cite{DGN} by a constant.

In \cite{DGN} and \cite{GP} a $C^2$ surface is called an H-minimal
surface if $H$ is identically zero.  In this paper, we make a
slightly different definition,
\begin{Def}  A $C^1$ graph $S$ over a domain $\Omega \subset \R^2$ is an {\bf H-minimal surface} if it
  satisfies equation \eqref{MSE} weakly.  More precisely, if
  $\nu_X=\pb \; X_1+\qb\; X_2$
  is the unit horizontal Gauss map of $S$, then $S$ satisfies
\[\int_\Omega \pb \phi_x+\qb \phi_y \;dx\;dy=0\]
for all $\phi \in C^\infty_0(\Omega)$.  
\end{Def}

In section \ref{glue}, we show that $C^2$ surfaces with $H=0$ are
H-minimal in this sense as well.  

For completeness, we also recall some of the results of \cite{GP}.  
\begin{Def}\label{seedcurve}
Let $S$ is a $C^2$ H-minimal graph and $\nux$ is its unit horizontal Gauss
map.  Thinking of $\nux$ as a vector field on $\R^2$, any integral
curve of $\nux$ is called a {\bf seed curve} of $S$.  We denote a seed
curve by $\gamma_z(s)$, i.e. $\gamma_z(0)=z,
\gamma_z'(s)=\nux(\gamma_z(s))$.  If a basepoint is understood, we
denote the curve by $\gamma(s)$.  We denote the integral curves of
$\nuxp$ by $\mathscr{L}_z(r)$ (or, simply $\mathscr{L}(r)$ if a
basepoint is understood).
\end{Def}

As mentioned in the introduction, in \cite{GP}, N. Garofalo and the
author show that, for $C^2$ H-minimal surfaces, $\mathscr{L}_z(r)$ are
straight lines in the plane and lift to horizontal lines in $\Hy$.
This yields an adapted parameterization of the plane:
\[ F(s,r) = (\gamma_1(s)+r\gamma_2'(s),\gamma_2-r\gamma_1'(s))\]
We recall that this parameterization ceases to be a local diffeomorphism
along the locus $r = \frac{1}{\kappa(s)}$ where $\kappa$ is the signed curvature of the seed curve $\gamma$ and
is given by 
\[\kappa(s) = \gamma''(s) \cdot \gamma'(s)^\perp \]

When lifted to $\Hy$, yields a parameterization of the
H-minimal surface as a ruled surface.

\begin{equation*}
S = \left ( \gamma_1(s)+r\gamma_2'(s),\gamma_2-r\gamma_1'(s),
  h_0(s)-\frac{r}{2}\gamma \cdot \gamma'(s) \right )
\end{equation*}

This is the content of theorem \ref{mrep} in the introduction.
Moreover, we can extend this parameterization from this patch of
surface to include all of the rules (i.e. allow $r \in (-\infty, \infty)$), which introduces characteristic
points at the locus given by:
\begin{equation}\label{cl} \Sigma(s,r)=h_0'(s)-r+\frac{1}{2} \gamma'(s) \cdot
\gamma(s)^\perp+\frac{r^2}{2} \kappa(s)=0 
\end{equation}

We recall that generically, this yields two branches of the
characteristic locus, one on one side of the locus
$r=\frac{1}{\kappa(s)}$ and the second on the other side of this
locus.  We refer the reader to \cite{GP}, section 7, for a more detailed
discussion of these features.

\section{Noncharacteristic $C^1$ H-minimal graphs}
 
In this section, we investigate $C^1$ H-minimal graphs.  We will
focus first on section of $C^1$ H-minimal graphs that do not have
characteristic points.  In this setting, we show that such graphs are
ruled surfaces as in the $C^2$ case.  At the end of the section, we
will address the question of characteristic points.  Throughout the
section, we will consider a surface defined by $(x,y,u(x,y))$ where
$u:\Omega \ra \R$ is a $C^1$ function.

\subsection{Weak directional derivatives}

At first, we will assume that the function $u$ defining the H-minimal
surface is at least $C^2$ and so the components of the unit horizontal
Gauss map are continuously differentiable.  Under this assumption,
we compute the directional derivative of $\pb$ in the direction of
$v=\left (1,-\frac{\pb}{\qb}\right )$ (the choice of this vector will become evident
in a moment).

\[D_v \pb=\nabla \pb \cdot \left (1, - \frac{\pb}{\qb}\right ) = \pb_x-\frac{\pb}{\qb}\pb_y
= \pb_x +\qb_y\]
The last equation is true because $\qb = \sqrt{1-\pb^2}$ and hence $\qb_y=-\frac{\pb}{\qb}\pb_y$.

Thus, we can interpret the integral equation
\[-\int_\Omega \pb \phi_x + \qb \phi_y\;dx\;dy =0\]
as a weak form of the equation \[\nabla \pb \cdot \left (1, -
  \frac{\pb}{\qb}\right )=0\]

In other words, if $S$ is an H-minimal surface, then $\pb$ is weakly
constant in the $\nuxp$ direction ($\nuxp$ and $\left
  (1,-\frac{\pb}{\qb}\right )$ point in the same direction).  We take
this as a definition:

\begin{Def}  The directional derivative of $\pb$ in the direction of
  $v=\left (1,-\frac{\pb}{\qb}\right)$ is weakly zero if
\[ -\int_\Omega \pb \phi_x + \qb \phi_y\;dx\;dy =0\]

In this case, we write $D_v \pb =0$
\end{Def}

\subsection{Rulings of $C^1$ H-minimal graphs}
Mimicing classical arguments, we have the following result.  
\begin{Lem}  Let $v =\left (1,-\frac{\qb}{\pb}\right )$ where
  $(\pb,\qb)$ is the horizontal Gauss map of an
  H-minimal graph, $S$, over a domain in $O \subset \R^2$ with $\pb,\qb \in W^{1,1}(O)$.  Assume that $\Omega\subset
  O$ be an open  domain with the following properties: 
\begin{enumerate}
\item The portion of $S$ over $\Omega$ has no
  characterstic points.
\item $v$ is continuous on $\Omega$ (i.e. $\qb
  \neq 0$).
\end{enumerate}
  Last, let $c_{x}(t)$ be an integral curve of $v$
  with $c(0)=x$ and let  $D^h_v \pb(x) = \frac{\pb(c_x(h))-\pb(c_x(0))}{h}$.
  Then, for $V \Subset \Omega$ and $h < dist(x,\partial
  \Omega)$, $D^h_v\pb(y) =0$ for a.e. $y \in V$.  
\end{Lem}

\pf  Let $w$ be a continuous vector field on $V$ and let $c^w_{x_0}(h)$ be an integral curve of $w$ passing through
the point $x_0$.  We note that $c^w_{x_0}(h)$, as a point set, coincides
with $\mathscr{L}_{x_0}(r)$ but is parametrized differently.  We may reparametrize $c^w_{x_0}$ so that
$(c^w_{x_0})'(s)=hw$.  Assuming briefly that $f$ is a smooth function,
we have
\begin{equation*}
\begin{split}
f(c^w_{x_0}(h))-f(c^w_{x_0}(0))&= \int_0^1 \nabla f(c^w_{x_0}(s)) \cdot
(c^w_{x_0})'(s)\; ds \\
&= \int_0^1 \nabla f(c^w_{x_0}(s)) \cdot
hw\; ds \\
&= h \int_0^1 D_w f(c^w_{x_0})(s) \; ds 
\end{split}
\end{equation*}
So, 
\[D^h_w f =\frac{f(c^w_{x_0}(h))-f(c^w_{x_0}(0))}{h} =\int_0^1
D_w f(c^w_{x_0}(s)) \; ds \]
Integrating appropriately, we have, for example, that
\[\int_V |D^h_w f|\; dx \le \int_V |D_w f| \; dx \]

Using standard mollification, we can smooth the function $\pb$
yielding a $C^\infty$ function $\pb_\epsilon$.  As we have restricted
to a set where $\pb$ is continuous (i.e. there are no characteristic
points), we know that $\pb_\epsilon \ra \pb$ uniformly as $\epsilon
\ra 0$.  Noting that $\overline{p}^2,\overline{q}^2,\overline{p}\overline{p}_y,\overline{q}\overline{q}_y \in L^1(\Omega)$, we have that $(\overline{p}^2+\overline{q}^2)_y =2\overline{p}\overline{p}_y+2\overline{q}\overline{q}_y=0$ in $L^1$.  Since $\frac{\overline{q}}{\overline{p}}$ is continuous on $\Omega$ we have $-\frac{\overline{q}}{\overline{p}}\overline{p}_y=\overline{q}_y$ in $L^1_{loc}(\Omega)$.  Thus, $D_v\overline{p}=\overline{p}_x+\overline{q}_y$ in $L^1_{loc}(\Omega)$.  Thus, by $H$-minimality, for $\phi \in C^\infty_0(\Omega)$, 
\[\int_\Omega D_v\overline{p} \phi = -\int_{\Omega}\overline{p}\phi_x+\overline{q}\phi_y =0\]
and we have $D_v\overline{p}=0$.

Next, let 
\[\pte = \frac{\pb_\epsilon}{\sqrt{\pb_\epsilon^2+\qb_\epsilon^2}},\]
\[\qte = \frac{\qb_\epsilon}{\sqrt{\pb_\epsilon^2+\qb_\epsilon^2}}\]
and 
\[ \tilde{v} = \left ( 1, -\frac{\pte}{\qte}\right)\]
As $\pb_\epsilon$ converges to $\pb$ uniformly on $V$, it follows by direct calculation that $\pte$ converges to $\pb$ uniformly on $V$ as well.  Similarly, $\qte$ converges to $\qb$ uniformly on $V$ and $\tilde{v}$ converges to $v$ uniformly on $V$.

Under this definition, we have that
\[\pte^2+\qte^2 = 1 \]
Differentiating with respect to $y$ and solving for $(\pte)_y$
we have 
\[(\pte)_y = -\frac{\qte}{\pte}(\qte)_y\]
So, 
\begin{equation*}
\begin{split}
D_{\tilde{v}} \pte &= (\pte)_x- (\pte)_y \frac{\pte}{\qte}
\\
&=(\pte)_x
+(\qte)_y\frac{\qte}{\pte}\frac{\pte}{\qte}\\
&=(\pte)_x+(\qte)_y
\end{split}
\end{equation*}

As $(\overline{p}_x)_\epsilon = (\overline{p}_\epsilon)_x,(\overline{q}_x)_\epsilon = (\overline{q}_\epsilon)_x$, it follows that $||(\overline{p}^2_\epsilon + \overline{q}_\epsilon^2)_x||_{L^1} \rightarrow 0$ as $\epsilon \rightarrow 0$.  So,
\[||(\tilde{p}_\epsilon)_x-\overline{p}_x||_{L^1(\Omega)} \le \sup \left | \frac{\overline{p}_\epsilon}{(\overline{p}^2_\epsilon + \overline{q}_\epsilon^2)^{\frac{3}{2}}}\right | || (\overline{p}^2_\epsilon + \overline{q}_\epsilon^2)_x||_{L^1(\Omega)} + \sup \left | \frac{1}{\sqrt{\overline{p}^2_\epsilon + \overline{q}_\epsilon^2}} \right | ||(\overline{p}_x)_\epsilon - \overline{p}_x||_{L^1(\Omega)} \rightarrow 0\]
Similarly, $||(\tilde{q}_\epsilon)_y-\overline{q}_y||_{L^1(\Omega)}\rightarrow 0$ and we conclude $||D_{\tilde{v}}\tilde{p}_\epsilon||_{L^1(\Omega)} \rightarrow 0$ as $\epsilon \rightarrow 0$.  Hence, there exists a function $C(\epsilon)\rightarrow 0$ as $\epsilon \rightarrow 0$ so that 
\[\int_\Omega |D_{\tilde{v}}\tilde{p}_\epsilon| \le C(\epsilon)\]

 Hence, there exists a function
$C(\epsilon)$ tending to zero as $\epsilon \ra 0$ so that 
\[ \int_V |D_{\tilde{v}} \pte| \le C(\epsilon)\]
So, applying the computation at the beginning of
the proof with $f=\pte$, we have 
\[\int_V |D^h_{\tilde{v}} \pte| \le C(\epsilon)\]

Thus, as $\epsilon \ra 0$, $D^h_{\tilde{v}} \pte \ra 0$ for almost every $x_0 \in
V$.  To complete the proof, we would like to have that $\lim_{\epsilon
  \ra 0} D^h_{\tilde{v}} \pte = D^h_v \pb$.  Assuming this for a
moment, this would imply that $D^h_v\pb=0$ almost everywhere as well.
Then, taking a countable dense sequence $h_n \ra 0$ and the countable 
intersection of full measure sets where $D^{h_n}_v \pb =0$, we have a
full measure subset of $V$, denoted by $V_0$, where
\[ D^{h_n}_v\pb(x_0)=0 \text{\;\; for all $n \in \mathbb{Z}_+,x_0\in
  V_0$}\]
By the continuity of $\pb$ on this region, this implies that 
\[D^h_v \pb(x)=0 \text{\;\; for all $h< dist(x,\partial \Omega), x \in V$}\]

Thus, we are left with verifying that 
\[
\lim_{\epsilon \ra 0} D^h_{\tilde{v}} \pte = D^h_v \pb
\]
First we note that $D^h_{\tilde{v}} \pte = D^h_{\tilde{v}} \pb +
o_\epsilon(1)$ and that the convergence is uniform as $\pte \ra \pb $
uniformly on $V$.  So, we merely need to verify that 
\begin{equation}\label{eqa}
\lim_{\epsilon \ra 0} D^h_{\tilde{v}} \pb = D^h_v \pb
\end{equation} 
To calculate the value of the limit, we will construct a sequences of
integral curves using the work in appendix \ref{app}.  Letting $\Omega = V$, $X=(\pb,\qb)$ and
$X_k=(\tilde{p}_{\epsilon_k},\tilde{q}_{\epsilon_k})$ for a sequence
$\epsilon_k \ra 0$, we apply the construction in appendix \ref{app} to
form appropriate integral curves for these vector fields.  Then lemma
\ref{l3} implies that $\pb(c^{X_k}_{x_0})(h) \ra \pb(c^X_{x_0}(h))$
and hence \eqref{eqa} is true.   
$\qed$

\begin{Rem}\label{r1}
We note that if we assume $v=\left( \frac{\qb}{\pb},-1\right )$ and that $S$
  is a noncharacteristic patch of surface where $v$ is continuous,
  essentially the same argument proves that 
\[D^h_v \qb(x) =0\] for $x \in V$ and $h < dist(x, \partial \Omega)$.
  We note that since $\pb^2+\qb^2=1$, if $D^h_v\qb=0$ implies that
  $D^h_v\pb=0$ as well.
\end{Rem}

\begin{Lem} If $S$, a $C^1$ H-minimal surface, is decomposed as $N \cup
  \bigcup_{i=1}^\infty K_i$, then on each $K_i$ with nontrivial interior and with $\pb,\qb \in W^{1,1}$, the integral curves of $\nuxp$ are straight lines. 
\end{Lem}

\pf  Let $\Omega_1 \subset K_i$ be the open set where $\qb \neq 0$ and let $\Omega_2 \subset K_i$ be the set where $\pb \neq 0$.  Then $\Omega_1 \cup \Omega_2 = K_i$.  Let $V_j$ be compactly contained in $\Omega_j$.  By the previous lemma,
since $D_v^h\pb=0$ on $V_1$, we have that for each
integral curve, $\mathscr{L}$ of $\nuxp$, $D_v^h\pb$ is zero on almost
every point of $\mathscr{L}$.  Thus, for these integral curves $\nuxp=(\qb,-\pb)$ is constant almost
everywhere along its own integral curves (which are the same as
the integral curves of $v$).  Thus, these integral curves are straight
lines except potentially on a set of measure zero.    However, the structure theorem of
Franchi, Serapioni and Serra Cassano says that $p$ is a continuous
function and hence, $\pb$ is discontinuous only at characteristic
points.  As the $K_i$ have no characteristic points, we see that $\pb$
is continuous on each $K_i$.  So, the integral curves are $C^1$ and
thus, since they are lines almost everywhere, they must simply be
straight lines. Similarly, the integral curves of $\nu_X^\perp$ are straight lines on $V_2$ as well using remark \ref{r1}.  Using a compact exhaustion of the $\Omega_i$, we see that the integral curves of $\nu_X^\perp$ on $K_i$ are straight lines.$\qed$

This is the same basic result we found in section 4 of \cite{GP} - the integral curves of $\nux^\perp$ are straight lines.  So, if we can construct a seed curve $\gamma$ as an
integral curve of the vector field $\nux$ and repeat the proof of Theorem 4.5 in \cite{GP} (this is theorem \ref{mrep} of the introduction), we have the same result for $C^1$ noncharacteristic H-minimal graphs:   

\begin{Pro}\label{repext}
Let $S \subset \Hy$ be a noncharacteristic patch of a $C^1$ H-minimal surface of the type
\[
S\ =\ \{(x,y,t)\in \mathbb H^1 \mid (x,y)\in \Om\ ,\ t = h(x,y)\}\ ,
\]
where $h : \Om \to \R$ is a $C^k$ function over a domain $\Om$ in the
$xy$-plane with $\pb,\qb \in W^{1,1}(\Omega)$.  Then, there exists a $C^1$ seed curve $\gamma$  so
that $S$ can be locally parameterized by
\begin{equation}
(s,r)\ \to\ (\gamma_1(s)+r\gamma_2'(s),\gamma_2(s)-r\gamma_1'(s),
h(s,r))\ ,
\end{equation}
where
\begin{equation}
h(s,r)\ =\ h_0(s)\ -\ \frac{r}{2} <\gamma(s), \gamma'(s)>\ .
\end{equation}
and $h_0(s)=h|_{\gamma(s)}$.  
\end{Pro}
We note that if we knew $\gamma \in C^2$, the argument used to prove
theorem 4.5 in \cite{GP} would extend completely to this case, showing
that a $C^1$ noncharacteristic graph is an H-minimal surface if and
only if it has such a representation for a neighborhood of each point
on the surface.  A priori, $\gamma$ is merely
$C^1$ and need not have any higher regularity.  However, we shall see
that if $\Omega$ is ``large in horizontal directions'' then $\gamma'$
is indeed $C^1$.  To make this precise, we need a definition.

\begin{Def}\label{hthick}  Suppose an open set $\Omega$ is parametrized by
\[F(s,r)= (\gamma_1(s)+r\gamma_2'(s),\gamma_2(s)-r\gamma_1'(s))\] 
where $\gamma$ is a seed curve.  Let \[d(s) = \min \{\sup \{r_1 |
r_1>0, F(s,r)|_{r \in (0,r_1)} \in
\Omega\}, \sup \{r_2 |
r_2>0, F(s,r)|_{r \in (-r_2,0)} \in
\Omega\}\}\] and let \[d(\Omega,\gamma)=inf \{d(s) | \gamma(s) \in \Omega\}\]
\end{Def}

\begin{Lem}\label{Lip} Fix $\epsilon >0$.  Let $S$ be a $C^1$
  noncharacteristic H-minimal graph defined over a planar domain
  $\Omega$ via $\gamma\in C^1$, a seed curve, and $h_0 \in C^1$, a height function,  for $S$.  If $d(\Omega,\gamma) >
  \epsilon$, then $\gamma'(s)$
  is locally Lipschitz.
\end{Lem}

\pf  We argue by contradiction.  Suppose $\gamma'(s)$ is not Lipschitz
at $s=s_0$.  Then, there exists a sequence $\{h_n\}$ tending to zero with 
\[\frac{|\gamma'(s_0+h_n)-\gamma'(s_0)|}{h_n} =\sqrt{2L_n}\]
with $L_n \ra \infty$ as $n\ra \infty$.  Now,
\[ |\gamma'(s_0+h_n)-\gamma'(s_0)|^2=
(\gamma'(s_0+h_n)-\gamma'(s_0))\cdot
(\gamma'(s_0+h_n)-\gamma'(s_0))=2-2\cos(\theta_n)\]
Where $\theta_n$ is the angle between $\gamma'(s_0)$ and
$\gamma'(s_0+h_n)$.  So, we must have that 

\[\frac{1-\cos(\theta_n)}{h_n^2} = L_n\]
Rearranging, we have
\begin{equation}\label{tn}
\cos(\theta_n) =1-h_n^2L_n
\end{equation}

\begin{figure}
\epsfig{file=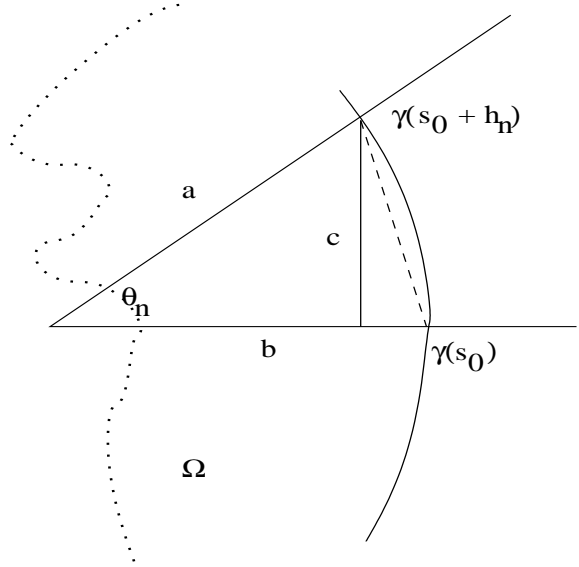,height=3in,width=3in}
\caption{Illustration for lemma \ref{Lip}}\label{lemfig}
\end{figure}

Recalling that, by proposition \ref{repext}, $\nux$ is constant along
integral curves of $\nuxp$, for $\nux$ to be well defined on $\Omega$,
no two integral curves of $\nuxp$ may cross inside $\Omega$.  Indeed,
it two such curves crossed, then infinitely many of them would cross and we would have conflicting values for $\nux$.

Next, we use this to gain an estimate on $L_n$.
Referring to figure \ref{lemfig}, we see that $\sin(\theta_n)
=\frac{c}{a} \le
\frac{|\gamma(s_0+h_n) - \gamma(s_0)|}{a}$.  Now, since
$d(\Omega)>\epsilon$, we have

\begin{equation}
\begin{split}
\epsilon &< d(\Omega) \le a\\
&\le \frac{|\gamma(s_0+h_n) - \gamma(s_0)|}{\sin(\theta_n)}\\
&= \frac{|\gamma(s_0+h_n) - \gamma(s_0)|}{\sqrt{1-(1-h_n^2L_n)^2}} \;
\; \text{(by equation \eqref{tn})}\\
\end{split}
\end{equation}

Thus, we have
\[\sqrt{1-(1-h_n^2L_n)^2} \le \frac{|\gamma(s_0+h_n) -
  \gamma(s_0)|}{\epsilon}\]
Or, after some algebraic simplification,
\[L_n \le \frac{|\gamma(s_0+h_n) -
  \gamma(s_0)|^2}{\epsilon^2h_n^2\left(1+\sqrt{1-\frac{|\gamma(s_0+h_n) -
  \gamma(s_0)|^2}{\epsilon^2}}\right)} \ra \frac{1}{2\epsilon^2}\]

But, by assumption, $L_n \ra \infty$ as $n \ra \infty$ so we reach a
contradiction.  $\qed$

\begin{Thm}\label{reg} Let $S$ be a $C^1$ H-minimal graph over a domain $\Omega$
of the xy-plane defined by $(\gamma,h_0)$ a seed curve/height function pair with $\gamma,h_0 \in C^1$.  If $d(\Omega,\gamma)>\epsilon>0$, then $\gamma$ is $C^2$.
\end{Thm}
\pf  By the previous lemma, we know that $\gamma'$ is Lipschitz and
hence $\gamma''(s)$ exists for almost every $s$.   Now, consider a
parameter value
$s_0$ and sequences $\{s_i^+\}$ and $\{s_i^-\}$ so that $s_i^+
\ra s_0$, $s_i^- \ra s_0$ and both $\gamma''(s_i^+)$ and
$\gamma''(s_i^-)$ exist for all $i$.  By the Lipshitz condition on
$\gamma'$, we can always find such sequences and we may also assume,
picking subsequences if neccesary, that $\lim_{i \ra \infty}
\gamma''(s_i^+)$ and $\lim_{i \ra \infty}
\gamma''(s_i^-)$ both exist.   Then, we claim that 
\[\lim_{i \ra \infty} \gamma''(s_i^+) = \lim_{i \ra \infty}
\gamma''(s_i^-)\]

To show this, we examine the Riemannian normal of the surface.  As the
surface is $C^1$, the normal must be continuous.  We will show that if
\[\lim_{i \ra \infty} \gamma''(s_i^+) \neq \lim_{i \ra \infty}
\gamma''(s_i^-)\]
then the normal cannot be continuous.

First, we a direct calculation using the representation of $S$ by \[\left (F(s,r),h_0(s)-\frac{r}{2} \gamma \cdot \gamma'(s) \right )\] yields that the vector
\[\eta(s,r) = \gamma_1'(s)\;X_1+\gamma_2'(s) \; X_2 + \beta(s,r) \; T\]
points in the direction of the Riemannian normal where 
\begin{equation}\label{beta}
\beta(s,r)=\frac{-1+r\kappa(s)}{h_0'(s) -r + \frac{1}{2}\gamma'\cdot
  \gamma^\perp + \frac{r^2}{2}\kappa(s)}
\end{equation}
We note that this computation is contained in section 4 of \cite{GP}.  
Now, by asumption, both $\gamma$ and $\gamma'$ are continuous.  To
argue by contradiction, we assume that 
\[\lim_{i \ra \infty} \gamma''(s_i^+)=l_1\]
and
\[\lim_{i \ra \infty} \gamma''(s_i^-)=l_2\]
where $l_1 \neq l_2$.  Let 
\begin{equation}\label{a2}
 \kappa_1=\lim_{i \ra \infty}\gamma''(s_i^+) \cdot
\gamma'(s_i^+)^\perp=l_1 \cdot \gamma'(s_0)^\perp
\end{equation}
and
\begin{equation}\label{bb3}
 \kappa_2=\lim_{i \ra \infty}\gamma''(s_i^-) \cdot
\gamma'(s_i^-)^\perp=l_2 \cdot \gamma'(s_0)^\perp
\end{equation}

Now, looking at $\eta(s,r)$ along the two sequences, we know that the
$X_1$ and $X_2$ components match as we tend towards $s_0$ as $\gamma'$
is continuous.  If $\eta$ is to be continuous, then $\beta$ must be
continuous as well, i.e.
\[ \frac{-1+r\kappa_1}{h_0'(s_0) -r + \frac{1}{2}\gamma'\cdot
  \gamma^\perp + \frac{r^2}{2}\kappa_1} = \frac{-1+r\kappa_2}{h_0'(s_0) -r + \frac{1}{2}\gamma'\cdot
  \gamma^\perp + \frac{r^2}{2}\kappa_2}\]

After simplifying algebraically, this yields:
\[(\kappa_1-\kappa_2)\left (rh_0'(s_0)+\frac{r}{2}\gamma' \cdot
  \gamma^\perp - \frac{r^2}{2}\right ) =0\]
As $r$ can vary, we see that $\kappa_1=\kappa_2$.  Since $\gamma$ is
parameterized by arclength, we have that $\gamma'' \cdot \gamma'=0$
where defined.  Combining this with equations \eqref{a2} and \eqref{bb3} we
reach a contradiction of the assumption that $l_1 \neq l_2$. So, we see that, where it is defined, $\gamma''$ coincides with a
continuous function.  Consider a point $s_0$ where, a priori, $\gamma'$
is not differentiable.  Then, in a neighborhood, $N$, of $s_0$, there
is a full measure subset $N_0$ so that if $s \in N_0$, $\gamma''(s)$
exists.  Then, as $\gamma''$ coincides with a continuous function where
is exists, we see that \[\lim_{s \in N_0,\; s \ra
  s_0}\frac{\gamma'(s_0)-\gamma'(s)}{s_0-s}\]
exists.  In other words, $s_0$ is a point of approximate
differentiability for $\gamma'$.  Since, by the previous lemma,
$\gamma'$ is Lipschitz, lemma 3.1.5 in \cite{Fed} implies that
$\gamma'$ is differentiable at $s_0$ as well. $\qed$  



\begin{Rem}  The previous theorem is a type of regularity result for
  H-minimal surfaces.  Recalling that $\nux = (\gamma_1'(s),\gamma_2'(s))$, the theorem says that the vector field $\nux$
  is continuously differentiable.  Therefore, the arguments from
  \cite{GP} used to prove
  theorem \ref{mrep} (this is theorem 4.5 in \cite{GP}) apply and the surface can be
  realized by 
\[\left (\gamma_1(s)+r\gamma_2'(s),\gamma_2(s)-r\gamma_1'(s),
h_0(s)-\frac{r}{2}\gamma \cdot \gamma'(s)\right)\]
Moreover, the smoothness of such a piece of H-minimal surface is completely
determined by the function $h_0(s)$. Given the structure theorem of
Franchi et al, if the surface is a perimeter minimizer, the function $h_0(s)$ must be at least $C^1$ on the
sets $K_i$. 
\end{Rem}

We end this section by summarizing the results:

\begin{Thm}  If $S$ is an open $C^1$ H-minimal graph over a domain
  $\Omega \subset \R^2$ with no
  characteristic points with unit horizontal Gauss map $\nux$ whose components are in $W^{1,1}(\Omega)$, then
  the integral curves of $\nux^\perp$ are straight lines and $S$ can
  be locally parameterized by
\begin{equation}\label{para}
(s,r)\ \to\ (\gamma_1(s)+r\gamma_2'(s),\gamma_2(s)-r\gamma_1'(s),
h(s,r))\ ,
\end{equation}
where
\begin{equation}\label{para2}
h(s,r)\ =\ h_0(s)\ -\ \frac{r}{2} <\gamma(s), \gamma'(s)>\ 
\end{equation}
and $\gamma$ is an integral curve of $\nux$.  Moreover, if there
exists $\epsilon >0$ so that $d(\Omega,\gamma)
>\epsilon >0$ then $h_0(s) \in C^1$ and $\gamma(s) \in C^2$.
\end{Thm}
This is theorem \ref{c1surf} in the introduction.

We note that, as $\gamma \in C^2$, all of the computations of section
4 of \cite{GP} are valid so long as they do not involve more than one
derivative of $h_0(s)$.  In particular, we have:
\begin{Pro}\label{surfext}  Let $S$ be a patch of $C^1$ H-minimal surface
  given by
 \begin{equation*}
(s,r)\ \to\ (\gamma_1(s)+r\gamma_2'(s),\gamma_2(s)-r\gamma_1'(s),
h(s,r))\ ,
\end{equation*}
where
\begin{equation*}
h(s,r)\ =\ h_0(s)\ -\ \frac{r}{2} <\gamma(s), \gamma'(s)>\ 
\end{equation*}
with $s \in (s_0,s_1), r \in (r_0,r_1)$.  Then, $S$ may be extended to
a surface $\tilde{S}$ by including all portions of the rules,
i.e. extending the above parameterization to $r \in
(-\infty,\infty)$.  In this case, the surface $\tilde{S}$ has
characteristic points at
\begin{equation}\label{cp}
h'_0(s)-r+\frac{1}{2}\gamma'(s) \cdot \gamma(s)^\perp +
\frac{r^2}{2}\kappa(s) =0
\end{equation}
\end{Pro}
\pf  The only new portion of this proposition is the identifiation of
the characteristic locus.  We note that by hypothesis, $d(\Omega) =
r_1-r_0 >0$ and so $\gamma \in C^2$.  As we assume the surface is
$C^1$, we must have that $h_0 \in C^1$ as well.  To verify the
position of the characteristic locus, we repeat the arguments found in
\cite{GP}, in particular the computations in the proof of theorem 4.6.  We review them here for completeness.  We
first compute tangent vectors to the surface at each point by taking
the $s$ and $r$ derivatives of the parameterization:

\begin{equation}
\begin{split}
\tau &= \frac{\partial}{\partial r}  (\gamma_1(s)+r\gamma_2'(s),\gamma_2(s)-r\gamma_1'(s),
h(s,r))\\
&= (\gamma_2'(s),-\gamma_1'(s), -\frac{1}{2} <\gamma(s),\gamma'(s)>)\\
&= \gamma_1'(s) \; X_1 + \gamma_2'(s) \; X_2 
\end{split}
\end{equation}

\begin{equation}
\begin{split}
\sigma &= \frac{\partial}{\partial s}  (\gamma_1(s)+r\gamma_2'(s),\gamma_2(s)-r\gamma_1'(s),
h(s,r))\\
&= (\gamma_1'(s) + r\gamma_2''(s),\gamma_2'(s)-r\gamma_1''(s), h_0'(s)-\frac{r}{2}-\frac{r}{2} <\gamma(s),\gamma''(s)>)\\
&= (\gamma_1'(s)+ r\gamma_2''(s)) \; X_1 +
(\gamma_2'(s)-r\gamma_1''(s)) \; X_2 +
(h_0'(s)-r\\&+\frac{1}{2}<\gamma'(s),\gamma(s)^\perp>+\frac{r^2}{2}
\kappa(s)) \; T 
\end{split}
\end{equation}
Taking the cross product of these vectors with respect to the
Riemannian structure, we have 
\begin{equation}
\begin{split}
\sigma \cross \tau &= \gamma_1'(s)B(s,r) \; X_1 + \gamma_2'(s)
B(s,r) \; X_2 + (-1 +r \kappa(s)) \; T
\end{split}
\end{equation}
where 
\[ B(s,r) = h_0'(s)-r +\frac{1}{2}<\gamma'(s),\gamma(s)^\perp>+\frac{r^2}{2}
\kappa(s)\]
As characteristic points arise when Riemannian normal, $\sigma \cross
\tau$, has only a $T$ component, we have the desired description of
characteristic points.  $\qed$

We often use the notation:
\[W_0(s)=h_0'(s) +\frac{1}{2}<\gamma'(s),\gamma(s)^\perp>\]

We end the section with an observation concerning the nature of the
characteristic locus along a single rule.  Equation \eqref{cl} shows
that, generically, each rule contains two characteristic points at 
\[ r=\frac{1}{\kappa(s)} \pm
\frac{\sqrt{1-2\kappa(s)W_0(s)}}{\kappa(s)}\]
one to each side of the point at $r= \frac{1}{\kappa(s)}$.  In the special case where $W_0(s)=\frac{1}{2
  \kappa(s)}$ we have a double characteristic point at $r=
\frac{1}{\kappa(s)}$.  

\begin{Lem}\label{rulecross}
Let  $S$ be a $C^1$ H-minimal graph parameterized by 
\[ \left (F(s,r), h_0(s)-\frac{r}{2} \gamma(s) \cdot \gamma'(s) \right), \; \; (s,r) \in \Omega \subset \R^2\]
Suppose $(s_0,r_0)$ and $(s_1,r_1)$ are points so that $F(s_0,r_0)=F(s_1,r_1)$.  Then 
$\left (F(s_0,r_0),h_0(s_0)-\frac{r_0}{2}\gamma(s_0)\cdot \gamma'(s_0)\right )$ is a characteristic point of $S$.
\end{Lem}

\pf  Since we assume that $S$ is a graph over the xy-plane, we must have that 
\[\left (F(s_0,r_0),h_0(s_0)-\frac{r_0}{2}\gamma(s_0)\cdot \gamma'(s_0)\right )=\left (F(s_1,r_1),h_0(s_1)-\frac{r_1}{2}\gamma(s_1)\cdot \gamma'(s_1)\right )\]
We recall that the unit horizontal Gauss map on $S$ is given by $\nux(s,r)
= (\gamma_1'(s),\gamma_2'(s))$ and that the unit horizontal Gauss map
is constant along any rule.  The vector 
\[ \eta(s,r)=\frac{\gamma_1'(s)}{\sqrt{1+\beta(s,r)^2}} \; X_1 + \frac{\gamma_2'(s)}{\sqrt{1+\beta(s,r)^2}} \; X_2 + \frac{\beta(s,r)}{\sqrt{1+\beta(s,r)^2}} \; T\]
where \[
\beta(s,r)=\frac{-1+r\kappa(s)}{W_0(s)-r+\frac{r^2}{2}\kappa(s)}\]
points in the same direction as the unit Riemannian normal to the surface
and is a completion of the unit horizontal Gauss map.  As the surface
is $C^1$, we must have that $\lim_{r \ra r_0} \eta(s_0,r)= \lim_{r \ra
  r_1} \eta(s_1,r)$.  Since we assume the two rules are not parallel
(they intersect), we have that $\gamma'(s_0) \neq \gamma'(s_1)$ and
hence, for these limits to be equal, we must have that \[\lim_{r \ra
  r_i} \beta(s_i,r) = \pm \infty\]

Examining the formula for the denominator of $\beta$ and equation
\eqref{cl}, we see that the intersection must be a characteristic
point.  $\qed$

\begin{Lem}\label{onecharpoint}  Let $S$ be a $C^1$ H-minimal graph
  parameterized by 
\[ \left (F(s,r), h_0(s)-\frac{r}{2} \gamma(s) \cdot \gamma'(s)
\right ), \; \; (s,r) \in \Omega \subset \R^2\]
Then along each rule, $\mathcal{L}_{\gamma(s)}(r)$, there is at most
one characteristic point.  
\end{Lem}

\pf As above, by equation \eqref{cl}, we see there are at most two
characteristic points along $\mathcal{L}_{\gamma(s)}(r)$.  Suppose
there are two characteristic points along a rule
$\mathcal{L}=\mathcal{L}_{\gamma(s_0)}(r)$, one to each side of
$r=\frac{1}{\kappa(s_0)}$.  We claim that, aribitrarily close to $r =
\frac{1}{\kappa(s_0)}$, $\mathcal{L}$ crosses another (nearby) rule.
To see this, we first left translate and rotate the Heiseneberg group
so that $\gamma(s_0)=0$ and $\gamma'(s_0)=(1,0)$ and reparametrize
$\gamma$ so that $s_0=0$.  From this
normalization, we have that $F(0,r)=(0,-r)$.  Consider a nearby
rule, $\mathcal{L}_{\gamma(s_1)}(r)$.  Then, direct calculation shows
that  
\[F\left(s_1,-\frac{\gamma_1(s_1)}{\gamma_2'(s_1)}\right)=F\left(0,\gamma_2(s_1)+\frac{\gamma_1(s_1)\gamma_1'(s_1)}{\gamma_2'(s_1)}\right)
\]
Taking that limit as $s_1 \ra 0$, we see that 
\[\gamma_2(s_1)+\frac{\gamma_1(s_1)\gamma_1'(s_1)}{\gamma_2'(s_1)} \ra
\frac{1}{\kappa(0)}\]
Thus, we make pick $s_1$ small enough so that
\[\mathcal{L}_{\gamma_{s_1}}(r) \cap \mathcal{L} \subset \mathcal{L}(r)|_{r \in \left ( \frac{1}{\kappa(0)}-\epsilon,
      \frac{1}{\kappa(0)}+\epsilon \right )}\]  By the previous lemma,
    we see that there must be a characteristic point at this
    intersection distinct from the two characteristic points assumed
    to be along $\mathcal{L}$.  This is a contradiction of \eqref{cl},
    which shows that there are at most two characteristic points.  Thus,
    along a rule contained in a piece of H-minimal graph, there is at
    most a single characteristic point.  $\qed$

\section{Continuous H-minimal surfaces}\label{glue}
Again taking our motivation from the theorem of Franchi, Serapioni and
Serra Cassano, we now investigate the possibility of gluing two pieces
of different of $C^1$ H-minimal surfaces together to form an new
H-minimal surface from their union.  In constrast to the classical
cases, we can, under certain restrictions, create piecewise $C^1$
surfaces that are globally merely continuous and yet satisfy the
H-minimal surface equation.  

We consider the problem, discussed in the introduction, of gluing
together two patches of $C^1$ 
H-minimal graphs so that the union satisfies the H-minimal surface
equation, at least weakly.

\noindent
{\bf Proof of theorem \ref{pwc1}:}  Assuming first that $S_1\cup S_2$ is H-minimal, we let 
\begin{equation*}
\nux = (\pb,\qb) = 
\begin{cases}
\nu_1 \text{\;\; on $\Omega_1$}\\
\nu_2 \text{\;\; on $\Omega_2$}
\end{cases}
\end{equation*}
Then,
\[ \int_{\Omega_1 \cup \Omega_2} \pb\phi_x+\qb\phi_y =0 \]
for a smooth compactly supported (on $\Omega_1 \cup \Omega_2$) test
function $\phi$.  Recall that by theorem \ref{reg}, we know that $\nux|_{\tilde{\Omega}_i}
\in C^1(\tilde{\Omega}_i)$.  First we compute
\begin{equation*}
\begin{split}
\int_{\Omega_i} \pb_i\phi_x+\qb_i \phi_y &= \int_{\Omega_i} (\pb_i
\phi)_x - \pb_{i,x} \phi + (\qb_i \phi)_y - \qb_{i,y}\phi \\
&= \int_{\Omega_i} (\pb_i
\phi)_x  + (\qb_i \phi)_y - \int_{\Omega_i}  \pb_{i,x} \phi + \qb_{i,y}\phi \\
&= \int_{\Omega_i} (\pb_i
\phi)_x  + (\qb_i \phi)_y \\
&= \int_{\partial \Omega_i} -\qb_i\phi \; dx + \pb_i \phi \; dy
\text{\;\; (by Green's theorem)} \\
&= \int_{C} -\qb_i\phi \; dx + \pb_i \phi \; dy \\
& = \int_C \phi (\nu_i \cdot \vec{n}_i) \; ds
\end{split}
\end{equation*}

The third equality holds because the surface over the interior of $\Omega_i$ satisfies
the minimal surface equation.  The second to last equality holds
because $\phi$ is compactly supported on $\Omega_1 \cup \Omega_2$ and
hence can only be nonzero on $C=\Omega_1 \cap \Omega_2$.  In the last equality, $\vec{n}_i$ denotes the inward pointing unit normal vector to $\partial{\Omega_i}$.

Applying this we have:

\begin{equation*}
\begin{split}
  \int_{\Omega_1 \cup \Omega_2} \pb\phi_x+\qb\phi_y &= \int_{\Omega_1}
  \pb_1\phi_x+\qb_1\phi_y + \int_{\Omega_2} \pb_2\phi_x+\qb_2\phi_y \\
&= \int_{\partial \Omega_1} \phi(\nu_1 \cdot \vec{n}_1) + \int_{\partial
  \Omega_2}  \phi(\nu_2 \cdot \vec{n}_2) \\
&= \int_{C} \phi(\nu_1 \cdot \vec{n}_1) + \int_{C}  \phi(\nu_2 \cdot
(-\vec{n}_1))\\
&= \int_C \phi((\nu_1-\nu_2) \cdot \vec{n}_1) \\
&=0
\end{split}
\end{equation*}
The second equality follows by the previous computation, where $\vec{n}_i$ are the inward pointing unit normal
  vectors.  The differentiability of the $\nu_i$ and the fact that the $S_i$ are H-minimal implies that  on the interior of $\Omega_i$, we have that
  $\pb_{i,x}+\qb_{i,y}=0$.  In the third equality, we observe that $\phi$ is zero on
  the boundary of $\Omega_1\cup\Omega_2$ and that on $C$,
  $\vec{n}_1=-\vec{n}_2$.  

Thus, we have that $(\nu_1-\nu_2)\cdot \vec{n}_1$ is weakly zero and
hence, $\nu_1-\nu_2$ must be tangent to $C$ almost everywhere.
Reversing the computation shows the sufficiency of this condition as
well.  $\qed$

We illustrate this with and example where we glue two different
H-minimal surfaces along a rule.  

\begin{eg}\label{eg1}  This theorem allows us to create many
  continuous H-minimal surfaces which are piecewise $C^1$.  We
  illustrate how to use this theorem by constructing a new H-minimal surface by gluing
  together the lower half of the plane $t=0$ with a portion of the
  surface $t=\frac{xy}{2}$ (see figure \ref{pw}).  To do this we
  define the following seed curve:
\begin{equation*}
\gamma(s)=
\begin{cases}
(1,s) \; \; 0 \le s < \infty \\
(\cos(s),\sin(s)) \; \; -\pi \le s < 0 \\
(-1, -s-\pi) \; \; -\infty < s < -\pi
\end{cases}
\end{equation*}
Note that, as a plane curve, $\gamma$ is two vertical lines glued to
the bottom half of a circle.  Now, we construct an H-minimal surface
from this seed curve.  With appropriate choices of $h_0(s)$, this yields
the parameterization:

\begin{equation*}
S:=
\begin{cases}
\left ( 1+r,s,-\frac{sr}{2}-\frac{s}{2}\right) \; \; \; 0 \le s < \infty, -1\le r < \infty\\
\left ( (1+r) \cos(s), (1+r) \sin(s), 0 \right ) \;\;\; -\pi < s < 0,
-1 \le r < \infty \\
\left ( -1-r, -s-\pi, -\frac{r(s+\pi)}{2}-\frac{s+\pi}{2} \right ) \;\;\; -\infty < s
\le -\pi, -1 \le r < \infty
\end{cases}
\end{equation*}

\begin{figure}
\epsfig{file=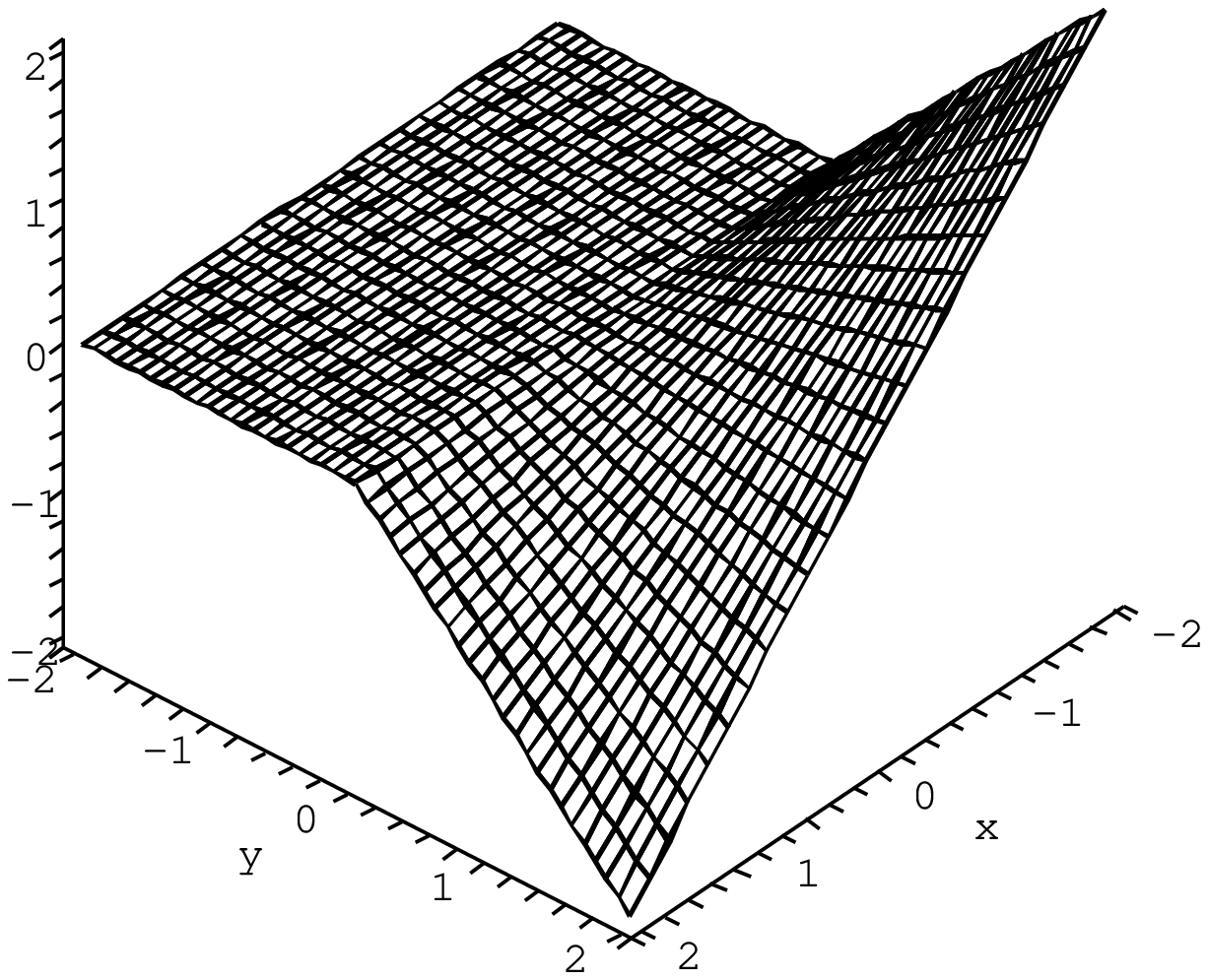,height=3in,width=3in}\label{pw}
\end{figure}

Calculating in xy-coordinates, we find the unit horizontal Gauss map
to be:

\begin{equation*}
\nux = (\pb,\qb)=
\begin{cases}
(0,\text{sgn}(x)) \; \; y \ge 0\\
\left ( -\frac{y}{\sqrt{x^2+y^2}},\frac{x}{\sqrt{x^2+y^2}} \right ) \;
\; y < 0
\end{cases}
\end{equation*}
Direct computation shows that, away from $y=0$, these are H-minimal
surfaces.  

We note that, using the notation of the theorem, $\nu_1=\nu_2$ along
the line $y=0$ and hence the hypotheses of the theorem are satisfied so
long as we pick $\Omega_1$ and $\Omega_2$ away from the characteristic
locus (for example, we might consider $\Omega_1=\{(x,y)|2\ge x\ge1,2\ge
y\ge0\}$ and
$\Omega_2=\{(x,y)| 2\ge x\ge 1, -2 \le y\le0\}$).
\end{eg}

\section{$C^\infty$ solutions to the Plateau Problem: persistent
  H-minimal surfaces}
In \cite{Pauls:minimal}, we showed that H-minimal graphs can arise as limits
of minimal surfaces in $(\mathbb{H}^1,g_\lambda)$.  In this section,
we examine those surfaces which are minimal for all values of $\lambda
\in [1,\infty)$.

\begin{Def}  A $C^2$ nonparametric graph is called a {\bf persistent
    H-minimal surface} if it is H-minimal and is  minimal in $(\mathbb{H}^1,g_\lambda)$
  for all $\lambda \in [1,\infty]$.
\end{Def}

As shown in \cite{Pauls:minimal} theorem 3.6 and 3.7 this implies that the graph is H-minimal
and, for any $C^2$ curve satisfying the bounded slope condition that such a surface spans, it is the solution to
the Plateau problem for this curve.  Thus, the persistent minimal
surfaces are a (small) class of smooth solutions to the Plateau
problem.  Moreover, as a consequence of Remark 1 in section 3 of
\cite{Pauls:minimal}, we have the following characterization of persistent
minimal surfaces.  

\begin{Thm}  An H-minimal surface $S = \{(x,y,h(x,y))\}$ is persistent if and
  only if $\Delta h =0$ where $\Delta$ is the usual Laplacian in
  $\R^2$.
\end{Thm}

In this section, we will use the representation formula from theorem
\ref{mrep} and some results from \cite{GP} to classify the
persistent H-minimal surfaces.  

\begin{Lem}  If an H-minimal graph $S$ is persistent then the
  signed curvature of its seed curve is constant.  
\end{Lem}

\pf First we assume that $S$ is a persistent H-minimal graph.  If $S$ is given by $(x,y,h(x,y))$ then $p=h_x-\frac{y}{2}$ and
$q=h_y+\frac{x}{2}$ and so $\Delta h=0$ is equivalent to $p_x+q_y=0$.
Using the notation from the previous section, we have $p=\alpha \pb$
and $q = \alpha \qb$ and so 
\begin{equation}
\begin{split}
p_x+q_y &= (\alpha \pb)_x + (\alpha \qb)_y \\
&= \nabla \alpha \cdot \nux + \alpha(\pb_x+\qb_y)\\
&= \nabla \alpha \cdot \nux \; \; \text{(since $S$ is H-minimal)}\\
&= 0 \; \; \text{(since we assume $S$ is persistent)}
\end{split}
\end{equation}

Thus, $\alpha$ is constant along the curves $F(s,t_0)$ and so we may
write $\alpha(s,r)=\alpha(r)$.  However, from theorem 7.1 in \cite{GP}, we
know that 
\[ \alpha(s,r) = \frac{\frac{\kappa(s)}{2}r^2-r+\alpha_0(s)}{1-\kappa(s)r}\]

Since $\alpha$ is constant along $F(s,0)$ this implies that
$\alpha_0(s) \equiv \alpha_0$ is constant and so, $\kappa(s)$ must
also be constant.
$\qed$

\begin{Thm}\label{pers}  The persistent H-minimal graphs fall into two
  categories up to isometric transformations of $(\mathbb{H}^1,\cc)$:
\begin{enumerate}
\item ($\kappa =0$) 
\begin{equation*}
\begin{split}
h(x,y)&=\frac{m}{1+m^2}(x-x_0)^2+\frac{m^2-1}{m^2+1}(x-x_0)(y-y_0) -
    \frac{m}{1+m^2}(y-y_0)^2+\\
   & \frac{a}{\sqrt{1+m^2}}(x-x_0)+\frac{am}{\sqrt{1+m^2}}(y-y_0)+b
\end{split}
\end{equation*}
for $m,a,b \in \R$.
\item ($\kappa \neq 0$) $S$, given in cylindrical coordinates is \[(\rho
  \cos(\theta),\rho \sin (\theta), a \theta + b)\] for $a,b \in \R$
\end{enumerate}
\end{Thm} 

\pf  By the previous lemma, we know that $\kappa$ must be constant for
$S$ to be persistent.  

\noindent
{\bf Case 1: $\kappa =0$} \\
If the curvature of the seed curve is zero, it is a line in the
plane.  By left translation, we may move the surface $S$ so that the
seed curve passes through the origin.  Thus, we may assume that
\[\gamma(s) = \left
  (\frac{s}{\sqrt{1+m^2}},\frac{ms}{\sqrt{1+m^2}}\right )\]
for some $m \in \R$.  Note that $\gamma(s)\cdot\gamma'(s) = s$  In this case the parameterization $F$ is simply a
linear transformation of the plane and we can write $(s,r)$ in terms
of $(x,y)$.  Indeed, we have
\begin{equation}
\begin{split}
s&=\frac{x}{\sqrt{1+m^2}}+\frac{my}{\sqrt{1+m^2}}\\
r&=\frac{mx}{\sqrt{1+m^2}}-\frac{y}{\sqrt{1+m^2}}
\end{split}
\end{equation}

Plugging this into the representation given in theorem \ref{c1surf}, we
get

\begin{equation}
\begin{split}
h(x,y) &= h_0(s)+\frac{1}{2}rs \\
&=  h_0 \left( \frac{x}{\sqrt{1+m^2}}+ \frac{my}{\sqrt{1+m^2}}
\right ) \\&\;\;\;\;+ \frac{1}{2}\left (
  \frac{mx}{\sqrt{1+m^2}}-\frac{y}{\sqrt{1+m^2}} \right )\left (
  \frac{x}{\sqrt{1+m^2}}+\frac{my}{\sqrt{1+m^2}}\right ) \\
&= h_0 \left( \frac{x}{\sqrt{1+m^2}}+ \frac{my}{\sqrt{1+m^2}}
\right ) +\frac{m}{1+m^2}x^2+\frac{m^2-1}{m^2+1}xy -\frac{m}{1+m^2}y^2
\end{split}
\end{equation}

Now, $S$ is persistent if and only if $\Delta h =0$,
\begin{equation}
\begin{split}
h_{xx}+h_{yy} &= \left (h_0'' \frac{1}{1+m^2}+2\frac{m}{1+m^2} \right )+
\left (h_0''
\frac{m^2}{1+m^2} - 2 \frac{m}{1+m^2} \right ) \\
&= h_0''\\
&=0
\end{split}
\end{equation}

Thus, surfaces in this case are persistent if and only if $h_0$ is
linear, i.e. $h_0 = cs+d$ for some $a,b \in \R$.  If we now left translate the resulting surfaces by a fixed
basepoint $(x_0,y_0,t_0)$, we have that the surface is given by:
\begin{equation}
\begin{split}
(x,y,h(x,y)) &= 
(x_0,y_0,t_0) \cdot \Bigg (x,y,c\left( \frac{x}{\sqrt{1+m^2}}+ \frac{my}{\sqrt{1+m^2}}
\right ) \\&+\frac{m}{1+m^2}x^2+\frac{m^2-1}{m^2+1}xy
-\frac{m}{1+m^2}y^2+d \Bigg )\\
&= \Bigg (x+x_0,y+y_0,t_0+ c\left( \frac{x}{\sqrt{1+m^2}}+ \frac{my}{\sqrt{1+m^2}}
\right ) \\&+\frac{m}{1+m^2}x^2+\frac{m^2-1}{m^2+1}xy
-\frac{m}{1+m^2}y^2+d + \frac{1}{2}x_0y - \frac{1}{2}y_0x \Bigg)
\end{split}
\end{equation}
Substituting $\overline{x} = x+x_0$ and $\overline{y}=y-y_0$ and
collecting terms yields the claim.  

\noindent
{\bf Case 2: $\kappa \neq 0$}\\

If $\kappa = c \neq 0$, then $\gamma(s)$ must be a circle and after a
suitable left translation, we may
write 
\[\gamma(s) = \left (\frac{1}{c}\cos(s),\frac{1}{c} \sin(s) \right )\]
Hence, $\gamma(s) \cdot \gamma'(s) =0$.  Moreover, the parameterization
$F$ yields 
\begin{equation}
\begin{split}
x &= \left ( r-\frac{1}{c} \right )\cos(s)\\
y &= \left ( r-\frac{1}{c} \right )\sin(s)\\
\text{And,} \\
s &= \arctan \left (\frac{y}{x} \right ) \\
r &= \sqrt{x^2+y^2}+\frac{1}{c}  
\end{split}
\end{equation}
Thus,
\[h = h_0(s) +\frac{r}{2}\gamma(s) \cdot \gamma'(s) = h_0\left (\arctan
  \left (  \frac{y}{x} \right ) \right ) \]
Computing the Laplacian of $h$ yields:
\[\Delta h = \frac{h_0''\left (\arctan
  \left (  \frac{y}{x} \right ) \right )}{x^2+y^2}\]
Thus, $\Delta h=0$ if and only if $h_0''(s)=0$ or that $h_0$ is
linear.  Thus, using cylindrical coordinates, $S$ is given by \[ (\rho
\cos(\theta), \rho \sin (\theta), a \theta +b)\] for $a,b \in \R$. $\qed$

We record the observation made above:

\begin{Cor}  Any $C^2$ closed curve satisfying the bounded slope condition lying on the surfaces identified
  in theorem \ref{pers} has a $C^\infty$ solution to the Plateau
  problem.  
\end{Cor}

\section{Obstructions to H-minimal spanning surfaces of high regularity}
Throughout the rest of this paper, we will be considering a smooth closed curve
\[c(\theta)=(c_1(\theta),c_2(\theta),c_3(\theta)) \subset \Hy\] with
the property that $c(\theta)$ is a graph over the projection of $c$ to
the xy-plane.  For ease of notation, we will denote this projection by
$\overline{c}(\theta) = (c_1(\theta),c_2(\theta),0)$.  When the
context is clear, we surpress the last coordinate of the
projection.  We will be considering H-minimal spanning surfaces for
these curves and moreover will consider only $C^1$ H-minimal spanning surfaces that
are ruled graphs.  To be precise, we make a definition:
\begin{Def}  A $C^1$ {\bf ruled H-minimal graph}, $S$, over a closed domain $\Omega \in \R^2$
  is an H-minimal graph with the property that every rule in $S$ that meets a the characteristic locus may be extended over the characteristic locus as a straight line.
\end{Def}

In other words, we will not consider gluings of
the type discussed in theorem \ref{pwc1}.  We
note that the work in \cite{GP} shows that $C^2$ H-minimal surfaces
satisfy this definition.

\begin{figure}\label{A0}
\epsfig{file=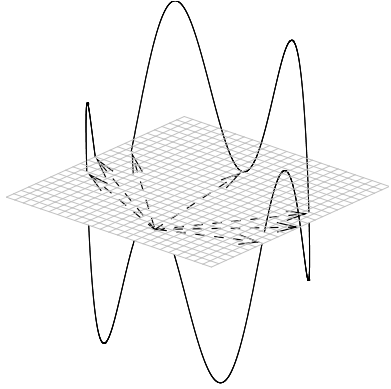}
\caption{The dotted arrows are the possible rules eminating from the point.}
\end{figure}

If $c$ lies on a $C^1$ ruled H-minimal graph then a geodesic line intersects each point on $c$ and, potentially, one or more other
points on $c$ (see figure \ref{A0}).  One easy way to determine the possible
geodesic lines which are allowable for a specfic point, $c(\theta_0)$, on the
curve is to left translate that point to the origin (recalling that left
translation preserves minimality).  At the origin, the horizontal
plane is the $xy$-plane and so, points which can be connected to
$c(\theta_0)$ by geodesic lines are those points on the left
translated curve which lie on the xy-plane.  Using the Campbell-Baker-Hausdorff
formula, one can calculate this set explicitly as:
\[ A(\theta_0)=\left \{\theta | c_3(\theta)-c_3(\theta_0)-\frac{1}{2}c_1(\theta)c_2(\theta_0)
  +\frac{1}{2}c_1(\theta_0)c_2(\theta) = 0\right \}\]

Note that $\theta_0 \in A(\theta_0)$.  In terms of building up a ruled
surface which spans $c$, the larger $A(\theta_0)$, the more
flexiblility one has in developing a surface.  On the other hand, if $A(\theta_0)$ contains only $\theta_0$ itself,
this places great restriction on the possibilities of smooth spanning
H-minimal surfaces.

Consider a $C^k$ closed curve $c:[t_0,t_1) \ra \Hy$ which is a graph
over a curve, $\overline{c}$, in the xy-plane.  Suppose $c$ is spanned by a ruled H-minimal
surface, $S$, which is a graph over a closed domain $\Omega$ in the xy-plane
so that $\partial \Omega = \overline{c}$.  Then the definition of $A$
above implicitly defines a function $\phi(t)$ for $t \in[t_0,t_1)$ via
the equation
\[c_3(\phi(t))-c_3(t)-\frac{1}{2}c_1(\phi(t))c_2(t)
  +\frac{1}{2}c_1(t)c_2(\phi(t)) = 0\]
As $c_i \in C^k$, $\phi$ is also $C^k$.  Moreover, we claim that such
a $\phi$ must be monotone.  To see this, suppose that $\phi$ is not
monotone and that there exist $t_0,t_1,t_2$ so that
$\phi(t_0)=\phi(t_2)$ and $\phi(t_1) \neq \phi(t_0)$.  Let
$\mathscr{L}_i$ be the rule connecting $c(t_i)$ to $c(\phi(t_i))$ and
let $\overline{\mathscr{L}}_i$ be the projection of $\mathscr{L}_i$ to
the xy-plane.  Then, the assumption on $\phi$ implies that
$\overline{\mathscr{L}}_1$ intersects either
$\overline{\mathscr{L}}_2$ or $\overline{\mathscr{L}}_0$.  Without
loss of generality, we will assume it intersects
$\overline{\mathscr{L}}_2$.  Further, $\overline{\mathscr{L}}_2$ must intersect the projection of
every rule connecting $c(t)$ to $c(\phi(t))$ for $t \in (t_0,t_1)$.
Such intersection points must be characteristic points of the surface
and so $\mathscr{L}_2$ would contain infinitely many characteristic
points in violation of lemma \ref{onecharpoint}.  

We record this observation:

\vspace{.25in}
\noindent
{\bf Existence Criterion:}  {\em Given a closed curve $c \in C^k$ which is a
graph over a curve in the xy-plane, if $c$ is spanned by a ruled
H-minimal graph then there exists a monotone $C^k$
function $\phi :S^1 \ra S^1$ with $\phi(\theta) \in A(\theta)$.  }

\vspace{.25in}

\begin{Def}  A point $c(\theta_0)$ is called {\bf Legendrian} if
  $c'(\theta_0) \in span \{X_1,X_2\}$.  We call a point
  {\bf isolated} if
  \[\{\theta_0\} = A(\theta_0)\]
\end{Def}

We record an immediate consequence of the definition.

\begin{Lem}\label{tang} If $c(\theta)$ is an isolated point and $c$
is spanned by a ruled H-minimal graph, then $c(\theta)$ is Legendrian and the rule passing through
$c(\theta)$ must be tangent to $c$.
\end{Lem}
\pf  If $c(\theta)$ is isolated then, by definition, it cannot be
connected to another point of $c$ via a rule.  As a consequence, we
note that the projection a rule through $c(\theta)$ to the xy-plane cannot be transverse to $\overline{c}$.
Indeed, if the projection were transverse, then it would
intersect another point on $\overline{c}$.  As $S$ is assumed to be a
graph, this rule would then be forced to intersect another point on
$c$.  Now, consider the rule through $c(\theta)$.  It is the limit of
rules connecting points near $c(\theta)$.  In other words, it is the
limit of secant lines and hence must be a tangent line to $c$ and $c(\theta)$.$\qed$

We next consider the relationship between two points on $c$ which are
connected by a rule on a spanning H-minimal surface.  

\begin{Lem}  Suppose $c(\theta_1)$ and $c(\theta_2)$ are connected by
  a rule, $\mathscr{L}$, of a ruled surface spanning $c$.  Then
  $c_3'(\theta_2)$ is proportional to the third coordinate of the parallel translation
  of $c'(\theta_1)$ along $\mathscr{L}$.  The proportionality constant
  is given as the derivative at the point $c(\theta_2)$ of the
  parametrization of the curve induced by the ruling around $c(\theta_2)$.
\end{Lem}

\pf  Without loss of generality, we may assume that $\theta_1=0$ and
$c(0)=0$ via a reparametrization of $c$ and a composition with left
translation in the Heisenberg group.  By hypothesis, $A(0)$ contains
the point $c(\theta_2)$ which, by abuse of notation, we will still
identify by the parameter value $\theta_2$ despite having
reparametrized the curve.  Note, that under this renormalization, the
rule $\mathscr{L}$ can be parametrized as 
\[\mathscr{L}(\tau) = (\tau c_1(\theta_2),\tau c_2(\theta_2), 0)\]
Moreover, the assumption that $c$ lies on a
ruled surface implies that there exists a mapping
$\phi:(-\epsilon,\epsilon) \ra (\theta_2-\delta,\theta_2+\delta)$
(with appropriately small $\epsilon$ and $\delta$) so
that $c(\phi(t)) \in A(c(t))$ for $t \in (-\epsilon,\epsilon)$ and
$c(t)$ is connected to $c(\phi(t))$ by a rule.  Thus, by the
definition of $A(\theta)$ we have
\begin{equation}
c_3(\phi(t))-c_3(t)-\frac{1}{2}c_1(\phi(t))c_2(t)+\frac{1}{2}c_1(t)c_2(\phi(t))
=0 
\end{equation}

Taking a derivative at $t=0$ and recalling that $c_1(0)=0=c_2(0)$ we
get
\begin{equation} \label{a1}
c_3'(\theta_2)\phi'(0) = c_3'(0)-\frac{1}{2}\overline{c}(\theta_2)
\cdot (c_2'(0),-c_1'(0))
\end{equation}

Next, we note that the vector field 
\[W=c_1'(0) \; X_1 + c_2'(0) \; X_2 + \left (c_3'(0)+\frac{\tau}{2}(\overline{c}(\theta_2)
\cdot (c_2'(0),-c_1'(0))\right ) \; T\]
is parallel along $\mathscr{L}$. Indeed, the tangent field to
$\mathscr{L}$ is given by \[V=c_1(\theta_2) \; X_1+ c_2(\theta_2) \;
X_2\] and a direct computation shows that $VW=0$.  Since the third
coordinate of $W$ at $\tau=1$ is proportional to the expression in equation (\ref{a1}),
we have the desired result.  $\qed$

\begin{Rem}  Geometrically, this says that the height function,
  relative to translation in the Heisenberg group, remains constant
  along the rules.  Thus, H-minimal surfaces are significantly more
  limited than ruled surfaces in Euclidean $\R^3$.  A comparable class
  of ruled surfaces in $\R^3$ would be those ruled surfaces that
  contain only rules parallel to the xy-plane.  
\end{Rem}

\subsection{Curves with isolated points}
We next turn to the problem of identifying curves that have ruled
H-minimal spanning surfaces.  We begin with an investigation of curves
that have isolated points.  
\begin{Lem}  Suppose $c$ is a $C^2$ curve in $\mathbb{H}^1$ and $c(\theta_0)$
  is an isolated point.  Then there is an open
  neighborhood, $I=(\theta_0-\epsilon,\theta_0+\delta)$, where $c|_I$ sits on a piece
  of a ruled surface. 
\end{Lem}

\pf  Note that, without loss of generality, by composing with a left
translation and reparametrizing the curve, we may assume that
$\theta_0=0$ and that $c(0)=0$.  We are attempting to parametrize a piece of the
curve for $t \in (-\epsilon,0]$ in terms of parameter values $t \in
[0,\delta)$ by associating a $\phi(t) \in (-\epsilon,0]$ to $t \in
[0,\delta)$.  So, we will construct a map $\phi:[0,\delta) \ra
(-\epsilon,0]$ with $\phi(0)=0$ so that $A(t)$ contains
$\phi(t)$.  By the definition of $A$, this implies that
\begin{equation}\label{aeq}
c_3(\phi(t))-c_3(t)-\frac{1}{2} c(\phi(t)) \cdot c(t)^\perp=0
\end{equation}

Differentiating with respect to $t$ solving for $\phi'(t)$, we get

\begin{equation}\label{phidot}
 \phi'(t) = \frac{c_3'(t)+\frac{1}{2} c(\phi(t)) \cdot
  c'(t)^\perp}{c_3'(\phi(t))-\frac{1}{2} c'(\phi(t)) \cdot
  c(t)^\perp} = \frac{c_3'(t)+\frac{1}{2} c(\phi(t)) \cdot
  c'(t)^\perp}{c_3'(\phi(t))+\frac{1}{2} c(t) \cdot c'(\phi(t))^\perp}
\end{equation}
Note that, at $t=0$, recalling that $c(0)=0$, we see that 
\[ \phi'(0)=  \frac{c_3''(0)}{\phi'(0) c_3''(0)} \implies \phi'(0)^2=1\]
Thus, for at least a
small time, $\phi(t)$, defined implcitly by ($\ref{aeq}$), exists and
hence, for some interval $I$, there exists a ruled surface spanning
$c|_I$.  $\qed$

\begin{Rem}  In the proof above, we see the obstruction - to be
  able to span a given curve with a ruled surface we must be able to
  find a function $\phi$ describing how to connect points on $c$ via
  rules that is monotone.
\end{Rem}

\begin{figure}[t]
\centering 
\mbox{\subfigure[$c_1(\theta)$]{\epsfig{figure=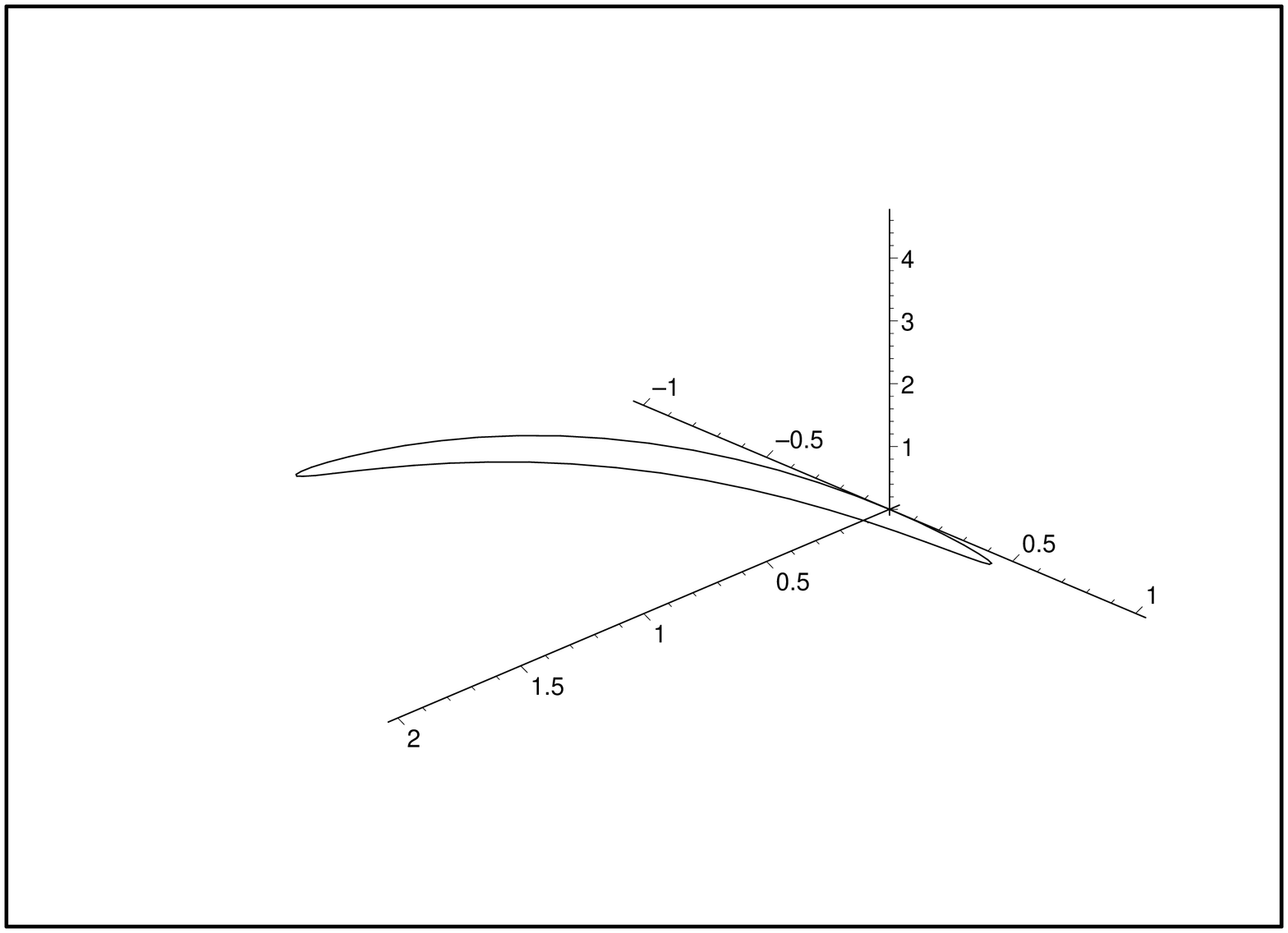,width=2.3in,height=2.3in,angle=-90}}\label{g1}\quad
\subfigure[$\theta \;\text{vs.}\; \theta_0$]{\epsfig{figure=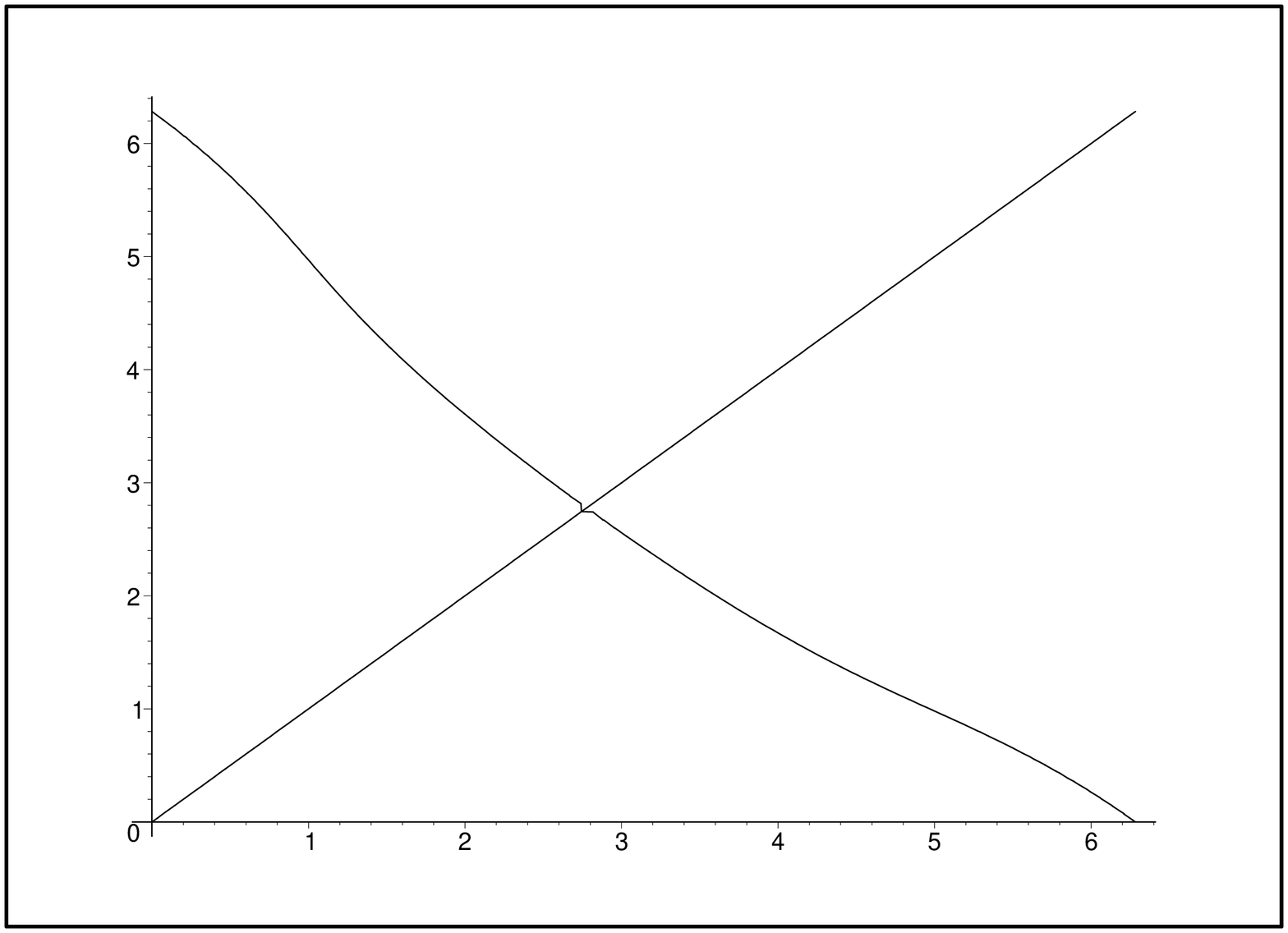,width=2.3in,height=2.3in,angle=-90}}}\label{g2}
\mbox{
\subfigure[Some rules of the spanning surface]{\epsfig{figure=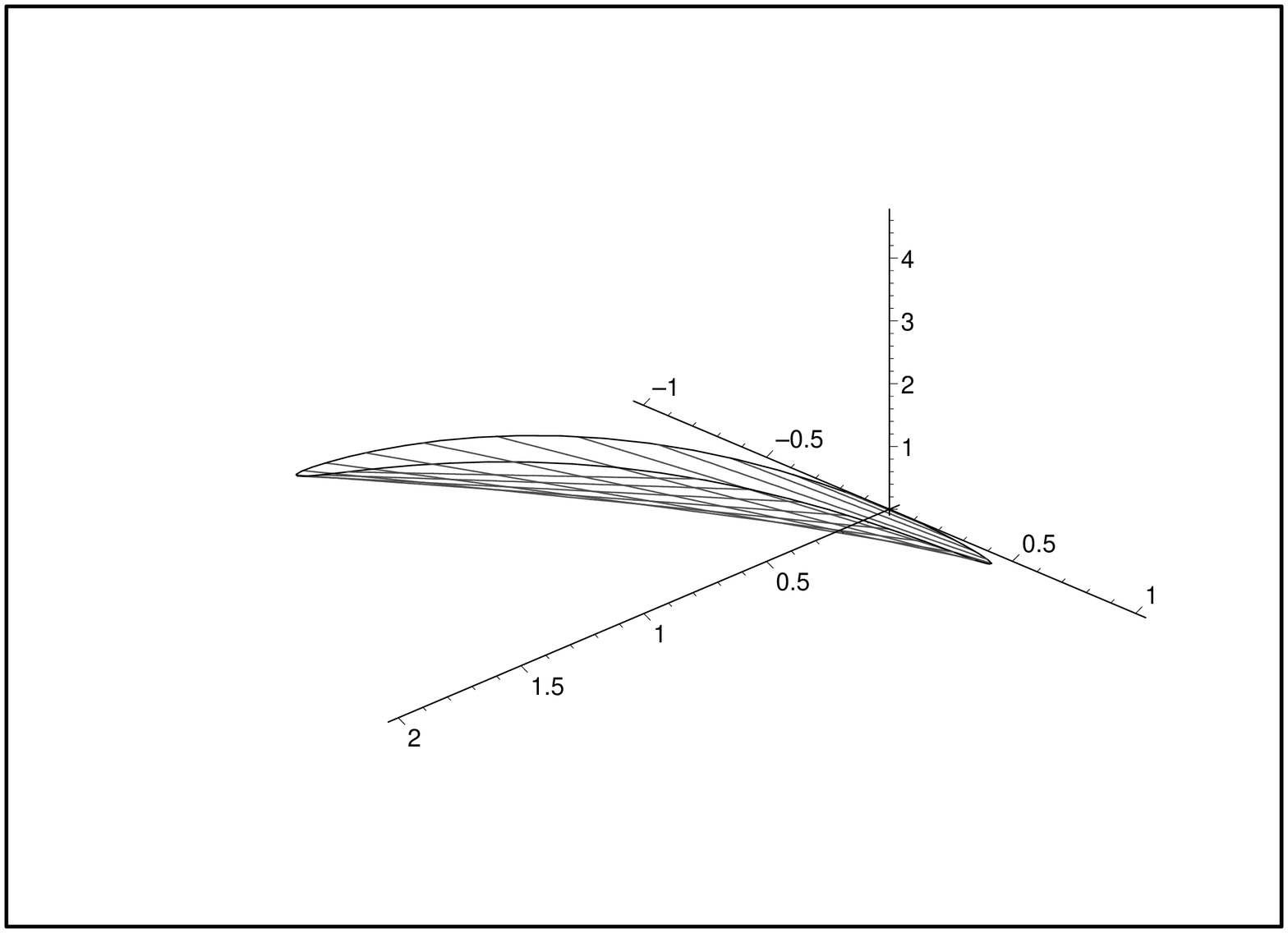,width=2.3in,,height=2.3in,angle=-90}}\label{g3}\quad
\subfigure[Projection to the xy-plane]{\epsfig{figure=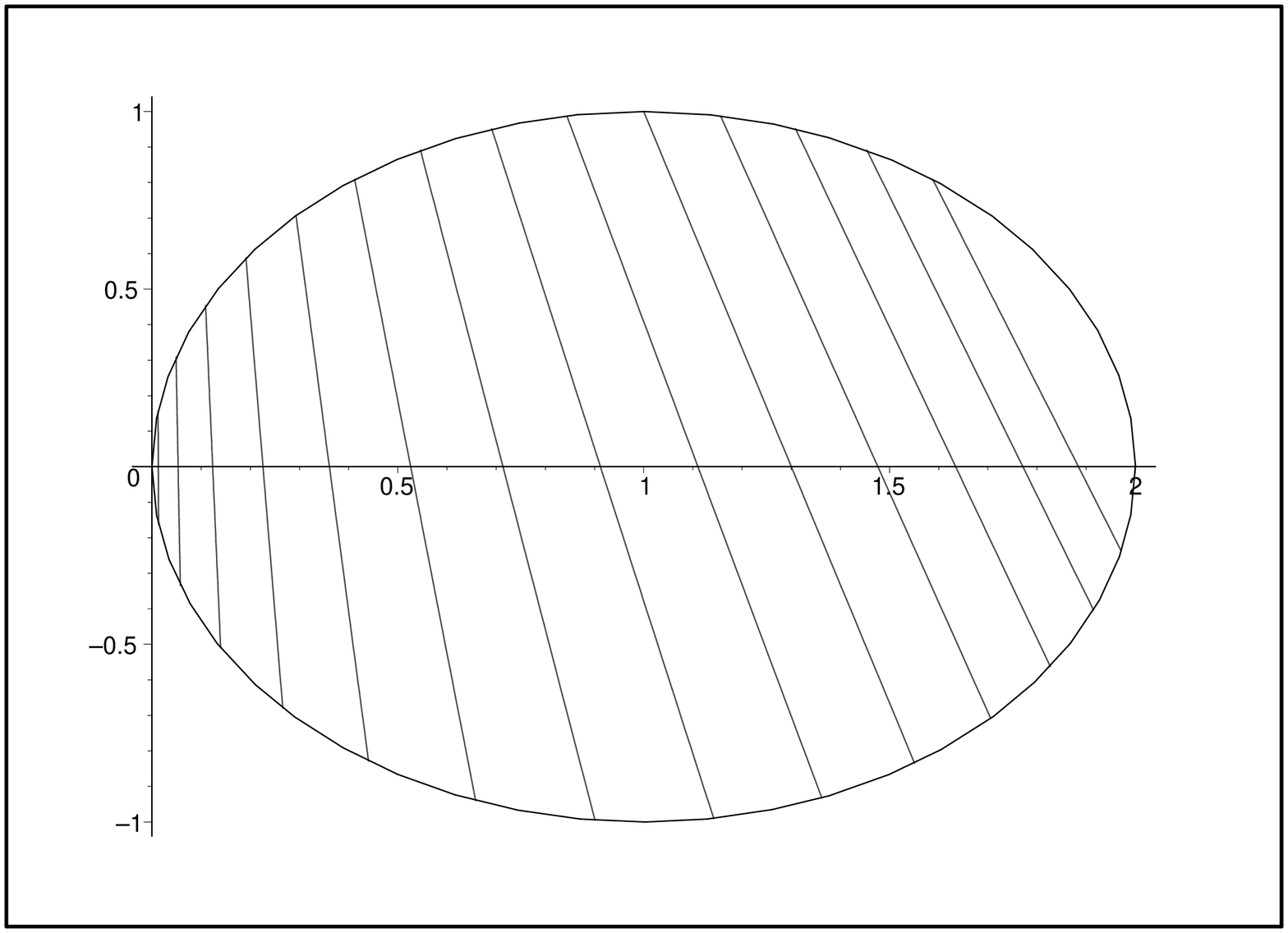,width=2.3in,,height=2.3in,angle=-90}}}\label{g4}
\caption{An example without an obstruction}\label{good}
\end{figure}

\begin{eg}  Consider the curve (see figure \ref{good}a)
\[ c_1(\theta) =
(1-\cos(\theta),\sin(\theta),2-2\cos(\theta)+\sin(\theta)-\sin(\theta)\cos(\theta))\]
for $\theta \in [0, 2\pi)$.  We quickly compute that 
\begin{equation*}
\begin{split}
A(\theta_0) = \{\theta | 2\cos(\theta_0)&+\frac{1}{2}
\sin(\theta_0)+\sin(\theta_0)\cos(\theta_0)-2\cos(\theta)\\ &+ \frac{1}{2}
\sin(\theta)-\sin(\theta)\cos(\theta)+\frac{1}{2}(\sin(\theta)\cos(\theta_0)
- \sin(\theta_0)\cos(\theta))\}
\end{split}
\end{equation*}
and note that for $\theta_0=0$,
\begin{equation*}
\begin{split}
A(0)&= \{\theta |
2-2\cos(\theta)+\sin(\theta)-\sin(\theta)\cos(\theta)=0\} \\
&= \{2n\pi\}
\end{split}
\end{equation*}
Thus $\theta_0=0$ is an isolated point.  Considering $\theta$ as a
function of $\theta_0$, we see in figure \ref{good}b that there is
another isolated point for $\theta_0$ slightly less than $\pi$.  We
will denote this value of $\theta_0$ by $\alpha$.  Observing figure
\ref{good}b, we see that for each
$\theta_0 = (0,\alpha)$ we can can connect $c_1(\theta_0)$ to the unique
point $c_1(\phi(\theta_0))$.  Figure \ref{good}c illustrates several of the
constructed rules connecting points on the curve and figure \ref{good}d
shows the projections of figure \ref{good}c to the xy-plane.

\end{eg}

\begin{figure}[t]
\epsfig{file=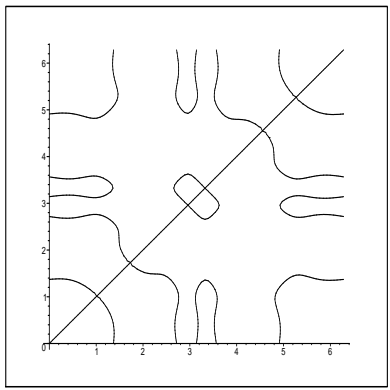,height=2in,width=2in,angle=-90}
\caption{A more complicated example}\label{b0}
\end{figure}

Of course, the example above is just about as well behaved as
possible.  However, the situation is often much more complicated.  For
example, if one considers the curve
\[ c(\theta) =
(1-\cos(\theta),\sin(\theta),\sin(4\sin(\theta)(1-\cos(\theta))))\]
figure \ref{b0} shows $\theta$ plotted as a function of $\theta_0$ (as
in figure \ref{good}b in the previous example).

\begin{figure}[t]
\centering 
\mbox{\subfigure[$c_2(\theta)$]{\epsfig{figure=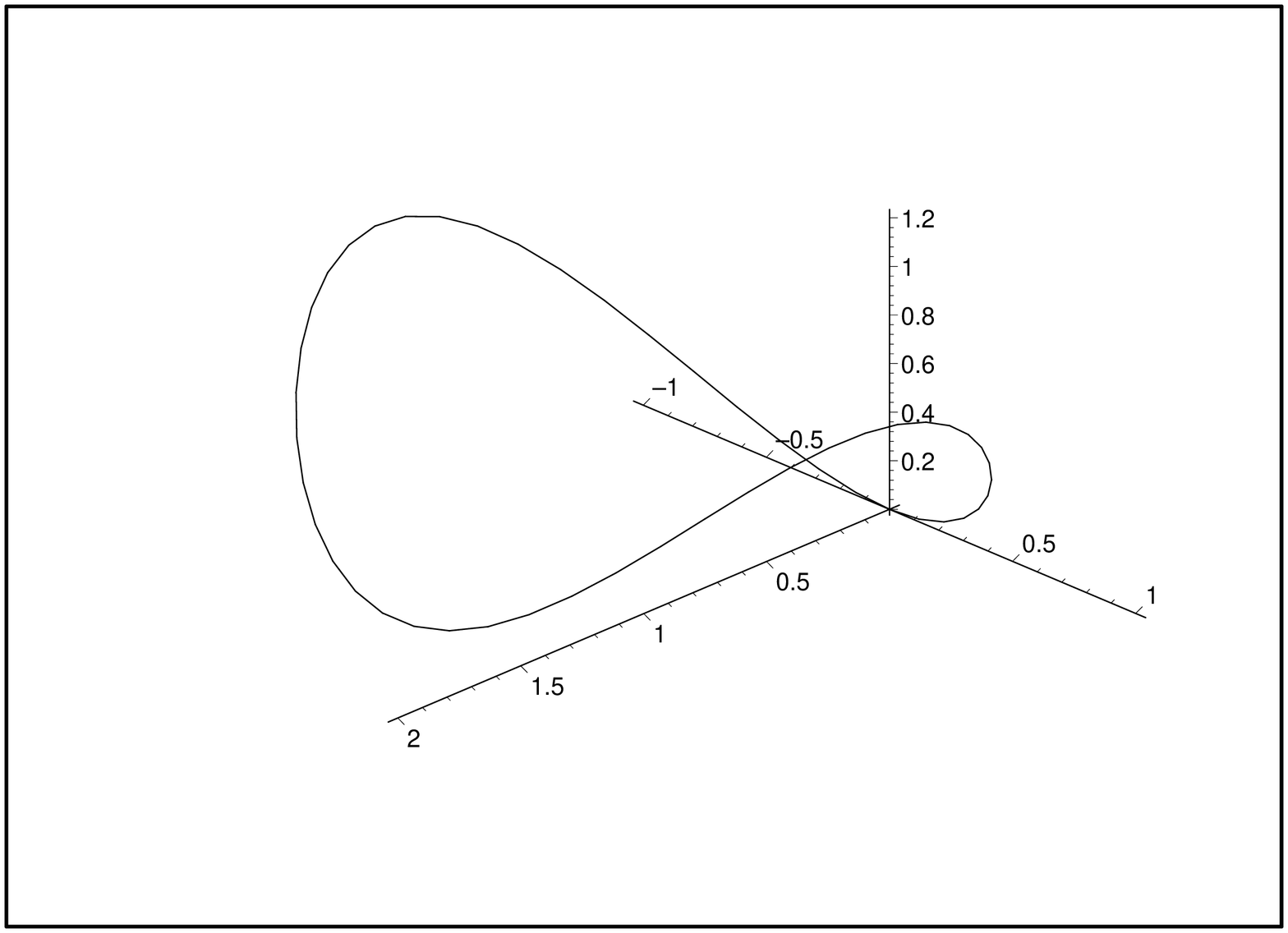,height=2.3in,width=2.3in,angle=-90}}\label{b1}\quad
\subfigure[$\theta \;\text{vs.}\; \theta_0$]{\epsfig{figure=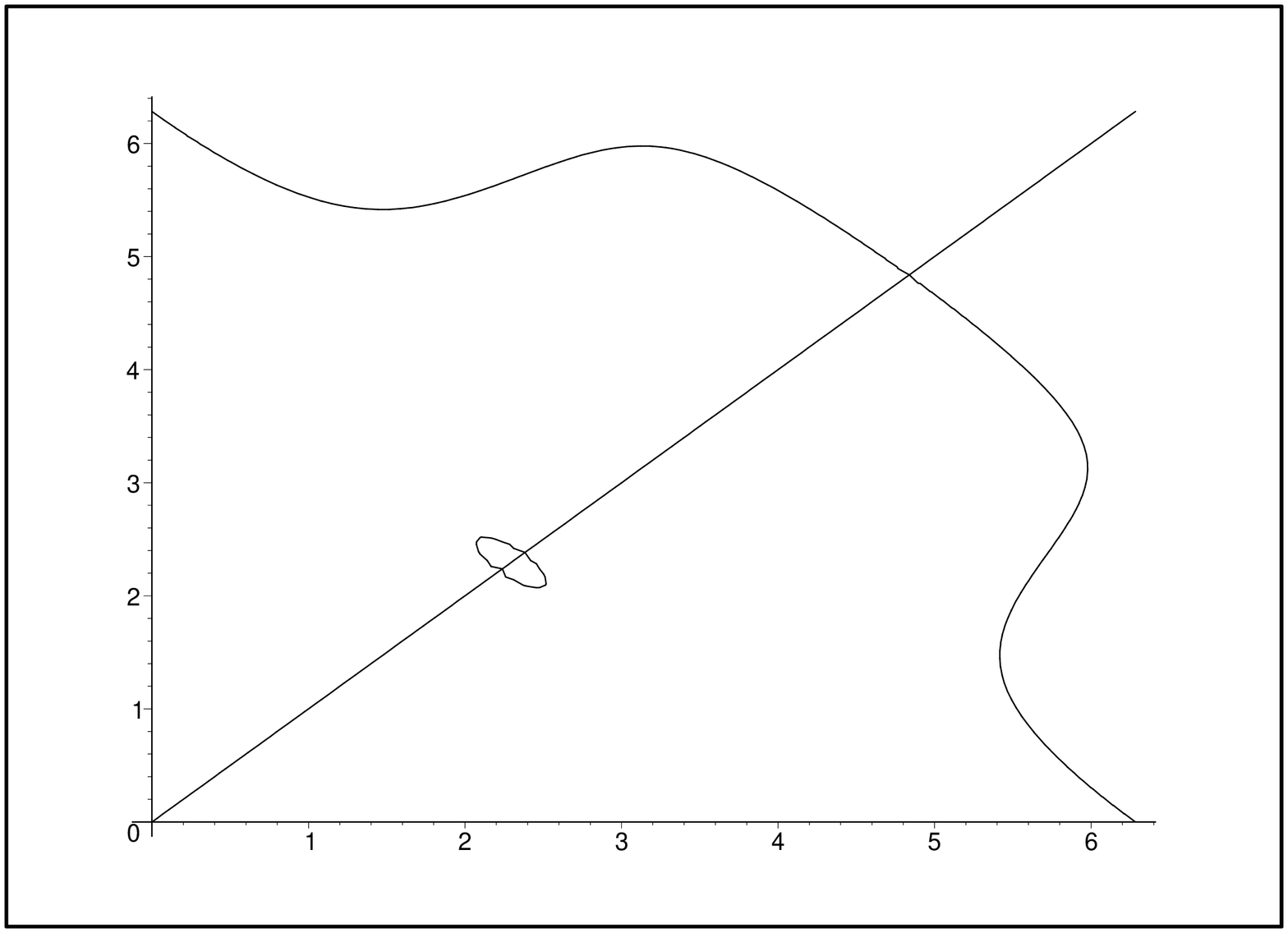,height=2.3in,width=2.3in,angle=-90}}}\label{b2}
\mbox{
\subfigure[Some rules of the spanning surface]{\epsfig{figure=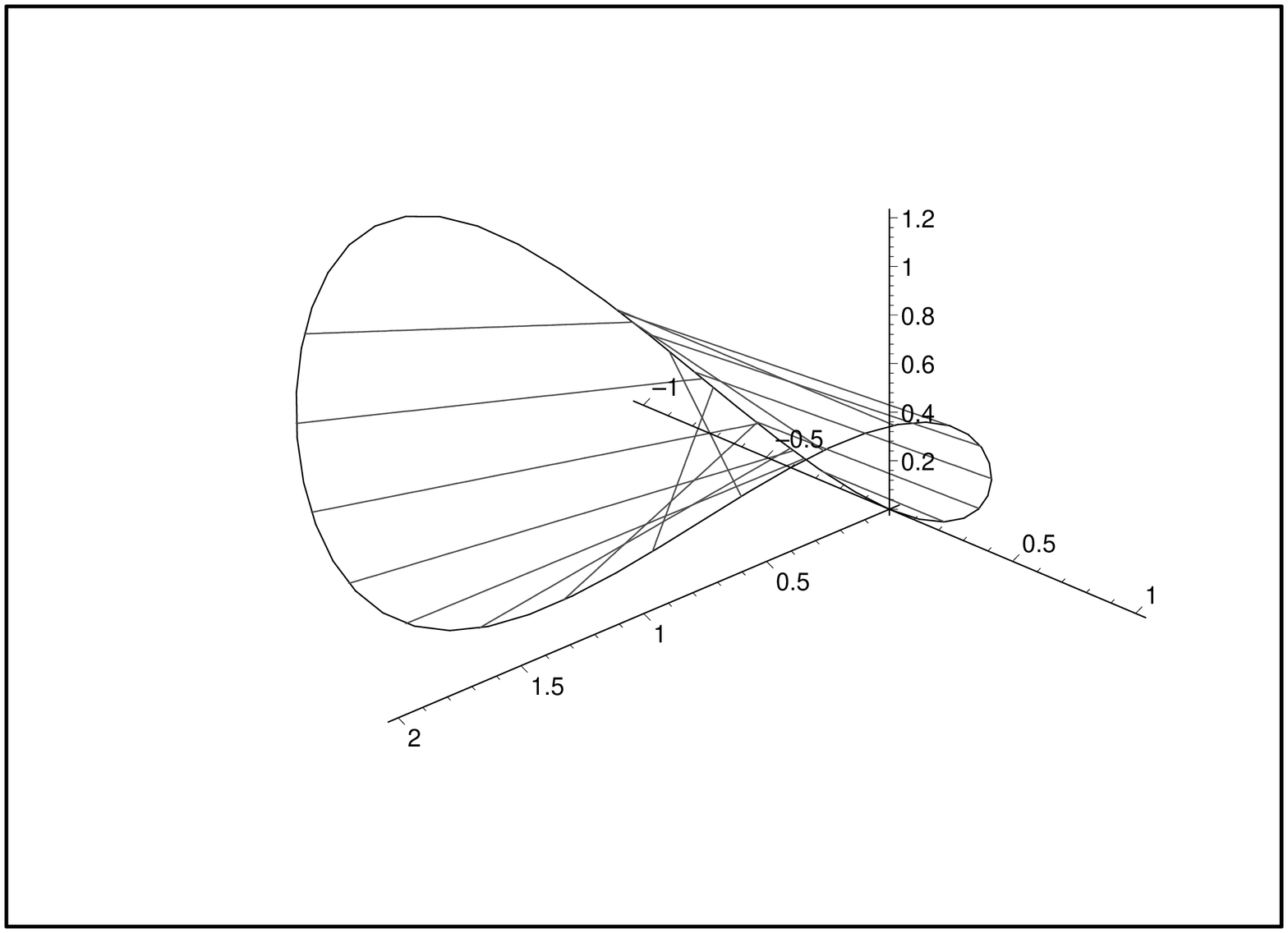,height=2.3in,width=2.3in,angle=-90}}\label{b3}\quad
\subfigure[Projection to the xy-plane]{\epsfig{figure=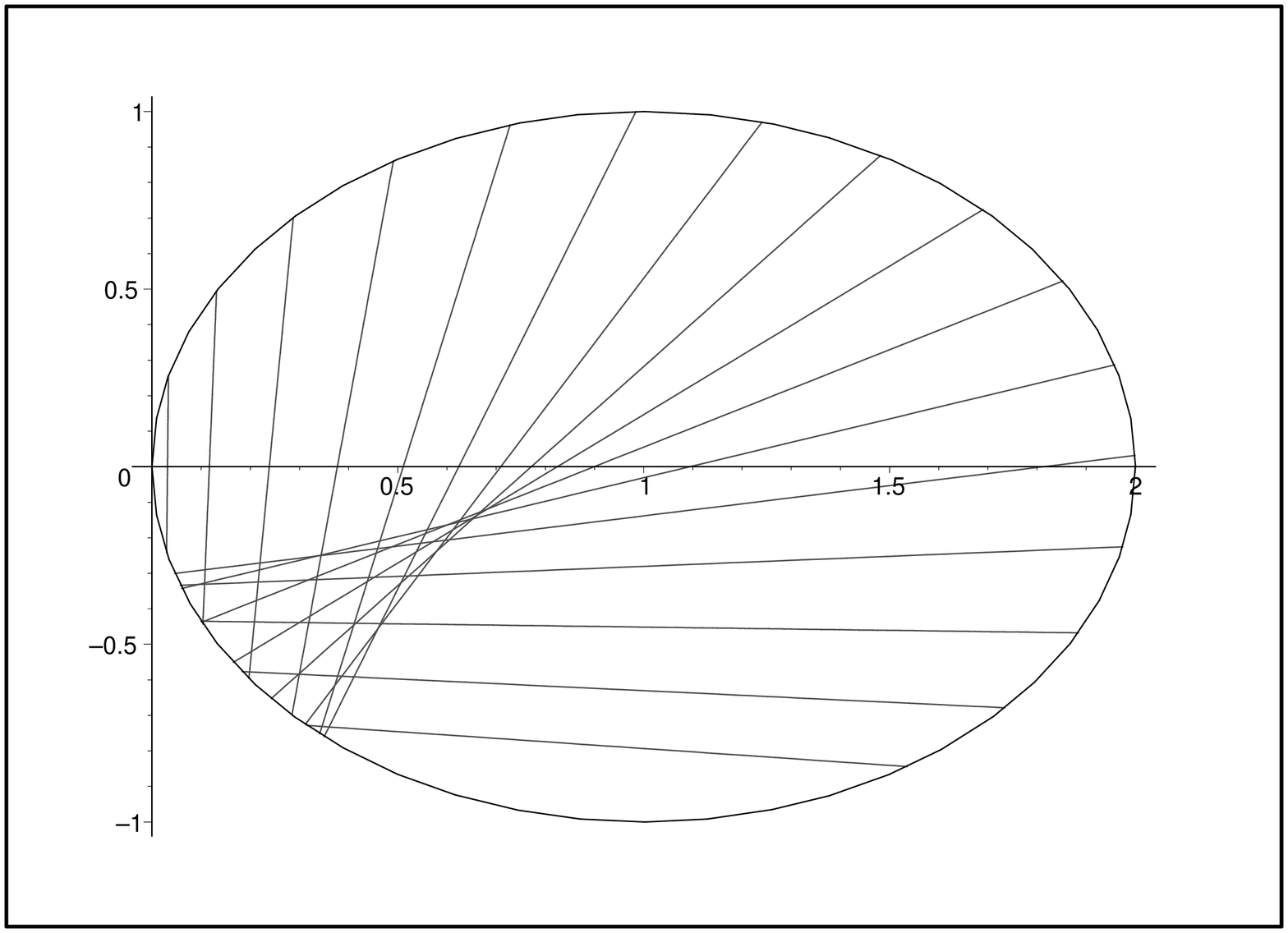,height=2.3in,width=2.3in,angle=-90}}}\label{b4}
\caption{An example with an obstruction}\label{bad}
\end{figure}

\begin{eg} In this example, we show an instance of the obstruction.
  In this case $\phi'(t)$ changes sign and we explore the effect of
  the sign change when building the ruled surface.  Consider the curve
  (see figure \ref{bad}a)

\[c_2(\theta) =
(1-\cos(\theta),\sin(\theta),1/5-1/5\cos(\theta)+\sin^2(\theta))\]

In this case, as in the first case, $\theta_0=0$ is an isolated point
and we can begin constructing a ruled surface.  Figure \ref{bad}b show
$\theta$ plotted as a function of $\theta_0$.  While tedious
computation can confirm this, the figure clearly shows that $\phi'$
changes sign at roughly $\frac{\pi}{2}$.  Figures \ref{bad}c and
\ref{bad}d illustrate how the construction fails - end endpoints of the
rules ``backtrack'' on the curve, creating a folded surface which, of
course, is no longer a graph.
\end{eg}

We end this discussion by noting the genericity of each of these
classes.  As the obstruction is defined by the monotonicity of $\phi$,
we note that strict monotonicity and non-mononicity are open conditions in
the $C^1$ topology by the implicit function theorem.  To make this
precise, we make the following definiton.  

\begin{Def}  A $C^2$ curve $c$ is {\bf generically nonmonotone} if
  there exists an $\epsilon > 0$ so that for
  any associated function $\phi$, there are parameter values $t_1,t_2$
  so that $\phi'(t_1) > \epsilon$ and $\phi'(t_2)< - \epsilon$.
\end{Def}
  
\begin{Pro}  Suppose for a given $c\in C^2$, there exists a $C^1$ ruled
  spanning H-minimal graph.  If the associated function $\phi$ is strictly
  monotone, then there exists an open neighborhood, $G$, of $C^2$ closed
  curves with respect to the $C^1$ topology containing $c$ so that any
  curve in $G$ has no obstruction to building a $C^1$ ruled H-minimal spanning surface. 

In addition, if $d$ is a $C^2$ generically nonmonotone curve, then there exists an open neighborhood of $d$ with
respect to the $C^1$ topology so that any curve in this neighborhood
cannot be spanned by a ruled minimal graph.
\end{Pro}

\pf  This is a consequence of formula \eqref{phidot}. For example, assume
that $\phi$ is strictly monotone decreasing:
\[ \phi'(t) =  \frac{c_3'(t)+\frac{1}{2} c(\phi(t)) \cdot
  c'(t)^\perp}{c_3'(\phi(t))+\frac{1}{2} c(t) \cdot c'(\phi(t))^\perp}
  <-a^2 <0 \]
Then, if we replace $c(t)$ with $\tilde{c}(t)=c(t)+\epsilon(t)$ where
  $|\epsilon(t)|< \delta$ and $|\epsilon'(t)| <\delta$, we have
\begin{equation*}
\begin{split}
 \phi_\epsilon'(t) &=  \frac{\tilde{c}_3'(t)+\frac{1}{2} \tilde{c}(\phi(t)) \cdot
  \tilde{c}'(t)^\perp}{\tilde{c}_3'(\phi(t))+\frac{1}{2} \tilde{c}(t)
  \cdot \tilde{c}'(\phi(t))^\perp}\\
& < -a^2 + o_\delta(1)
\end{split}
\end{equation*}

Thus, for $\delta$ sufficiently small (i.e. $\tilde{c}$ sufficiently
close to $c$ in the $C^1$ topology), $\phi_\epsilon'$ is strictly
monotone decreasing.  A similar argument shows the same genericity result for
curves where $\phi$ is generically nonmonotone.  $\qed$

\subsection{Totally non-Legendrian curves}
\begin{Thm}\label{nonL}  Suppose $c$ is a $C^1$ curve with no Legendrian points
  which is contained in an open $C^1$ ruled H-minimal graph, $S$.  Then there
  exists an interval, $I$, so that $c(I)$ is contained in a plane. 
\end{Thm}
\pf  We first record some easy observations:
\begin{itemize}
\item $S$ cannot have a characteristic point at any point of $c$.  If
  $c(t_0)$ were a characteristic point, then any smooth curve through
  $c(t_0)$ would be tangent to $\mathcal{H}_{c(t_0)}$, including $c$
  itself.
\item Consider a point $c(t_0)$ and let $\gamma$ be a seed curve
  through $c(t_0)$.  Use theorem \ref{c1surf}, to parameterize a neighborhood,
  $N$, of $c(t_0)$ that includes $c(t)$ for $t \in J$ where $J$ is an
  appropriate interval containing $t_0$.  By the nonLegendrian
  assumption and continuity of the normal vector, we may assume that there are no characteristic points in
  $N$ and hence, by theorem \ref{reg}, $\gamma \in C^2$.  Using the
  parametrization given by theorem \ref{c1surf}, there exist functions $s(t)$ and $r(t)$ so that $c(t) \cap
  N$ is parametrized by 
\[\left (\gamma_1(s(t))+r(t)\gamma_2'(s(t)),
  \gamma_2(s(t))-r(t)\gamma_1'(s(t)), h_0(s(t))-\frac{r(t)}{2}\gamma
  \cdot \gamma'(s(t)) \right)\] 
\item A rule through $c(t_0)$ is transverse to $c(t)$.  Indeed, if the
  rule were tangent then, by definition, $c(t_0)$ is
  Legendrian.
\item  For every $t \in J$, $\gamma(s(t))$ is twice
  differentiable and, applying formula \eqref{cp} at these points
  determines the charactersitic points along the rules passing through
  those points.
\end{itemize}

Let $\mathscr{L}_t(r)$ be the rule
through $c(t)$ and let $\overline{\mathscr{L}}_t(r)$ be the projection
of the rule to the xy-plane.   

\vspace{.25in}

\noindent
{\bf Claim:} {\em There exists two rules of $S$ that intersect in the
  interior of the portion of $S$ bounded by $c$}

\vspace{.25in}

\begin{figure}
\epsfig{file=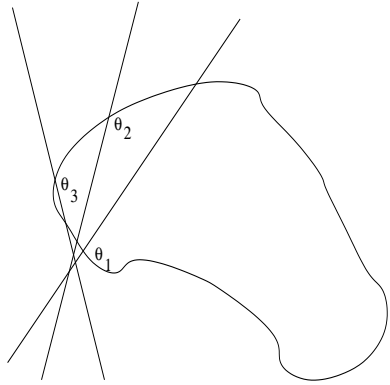,height=2in,width=2in}
\caption{Hueristic for picking points in the proof of Theorem \ref{nonL}} \label{cutfig}
\end{figure}

To show the claim, we assume there are no such rules.   First pick a parameter value $\theta_1$ and let
$\mathscr{L}_{\theta_1}(r)$ be the rule through $c(\theta_1)$.  Under
the assumption that the rule does not intersect any other rules inside
$c$, it must intersect another point on $c$, dividing $c$ (and the
surface) into two parts.  Next, pick a parameter value, $\theta_2$,
so that $c(\theta_2)$ is on the ``left hand side'' of the cut (see
figure \ref{cutfig}).  The rule, $\mathscr{L}_{\theta_2}(r)$ again must
cut the remaining portion into two parts.  We continue this iterative process, picking a sequence of parameter values $\{\theta_i\}$.
By construction, this sequence must converge to a value
$\theta_\infty$.  Moreover, if the rule $\mathscr{L}_{\theta_\infty}(r)$
does not intersect any of the $\{\mathscr{L}_{\theta_i}(r)\}$, it must
be tangent to $c$ at $\theta_\infty$.  This implies that
$c(\theta_\infty)$ is a Legendrian point, violating the hypothesis.  

Now, by the claim, we can pick $t_1,t_2$ so that $\overline{\mathscr{L}}_{t_1}(r) \cap
\overline{\mathscr{L}}_{t_2}(r) \neq \emptyset$.  Then the projection of
these two rules must not be parallel and hence, as the projections are
lines, they must intersect in a single point, $\{\overline{x}\}$.  As
$S$ is a graph over the xy-plane, we see that
$\mathscr{L}_{t_1}(r)\cap \mathscr{L}_{t_2}(r) = \{x\}$ where $x$ is
the point on $S$ over $\overline{x}$.  By lemma \ref{rulecross}, $x$ must be a
characteristic point of $S$.  A consequence of this observation is
that 
\[ \bigcap_{t \in J} \mathscr{L}_t = \{x\}\]

Suppose that this claim is not true, i.e. that there exists $t_0$ so
that $\mathscr{L}_{t_1}(r)\cap \mathscr{L}_{t_0}(r) = \{x'\}\neq
\{x\}$.  Then along $\mathscr{L}_{t_1}$ there must be 2 characteristic
points.  By lemma \ref{onecharpoint}, as $S$ is a graph over the xy-plane,
this cannot happen.  Using a left translation, we
may assume that $\{x\}$ is the origin.  Let $I$ be the interval
between $t_0$ and $t_1$.

Denote the portion of
$S$ bounded by $c(I)$, $\mathscr{L}_{t_0}(r)$ and
$\mathscr{L}_{t_1}(r)$  by $S_0$.  We finish the proof by showing that
$S_0$ is a portion of a plane.
Since the origin is a characteristic point, $S$ must be tangent to the
xy-plane at the origin.  As each rule is a horizontal straight line and every
rule in $S_0$ passes through the origin, we have that every such rule lies in the
xy-plane.  Thus, $S_0$ is a piece of the xy-plane and so $c(I)$ is planar.
$\qed$

\begin{Cor}  If $c$ is a smooth curve with no Legendrian points and no
  portion of 
  $c$ is contained in a plane then $c$ cannot be spanned by a $C^1$ ruled H-minimal surface.
\end{Cor} 

In particular, the solution to the Plateau Problem for such a curve
cannot be $C^2$ and, if it $C^1$, cannot be a ruled surface.  The best
result in this case would be a $C^1$ H-minimal surface composed of
ruled $C^1$ H-minimal patches glued along a mutual
characteristic locus so that the rules do not extend over the
characteristic locus.  One aspect of this type of surface is that the
horizontal Gauss map will necessarily have a discontinuity over the
characteristic locus which cannot be resolved by picking a horizontal
orientation.    

\begin{eg}
Let \[c(\theta)=(1-\cos(\theta),\sin(\theta),f(\theta))\]  Then,
\[c'(\theta)= \sin(\theta) \; X_1 + \cos(\theta) \; X_2 +\left (
f'(\theta)-\frac{\cos(\theta)}{2}+\frac{1}{2} \right ) \; T \]

Thus, for any periodic $f$ so that 
\[\left |f'(\theta)-\frac{\cos(\theta)}{2}+\frac{1}{2}\right | > \epsilon > 0\]
for some fixed $\epsilon$, we have an example that is totally non-Legendrian.

An explicit example is given by 
\[f(\theta)=\frac{1}{2} \sin(\theta)+ \frac{1}{8}\sin^2(\theta)\]
It is easy to show that no portion of this curve is planar.  
\end{eg}

\appendix
\section{Integral curves of continuous vector fields}\label{app}

In this appendix, we will review the existence of integral curves for
continuous vector fields and prove some results needed in the main
body of the paper.  The results here are consequences of Picard's
standard iterative construction of solutions to first order ordinary
differential equations (see, for example, \cite{Hartman}).  Our only
modification is to restrict our view to merely continuous vector
fields (as opposed to Lipschitz continuous ones). We make the
following standing assumptions:
\begin{enumerate}
\item  Let $X$ be a vector field defined on a compact
domain $\Omega$.
\item Let $\{X_k\}$ be a sequence of $C^\infty$ vector
fields, defined on $\Omega$, which converge uniformly to $X$ on
$\Omega$.
\end{enumerate}  

Let $\{M_k\}$ be a sequence constants tending to zero so
that $|X_k-X| \le M_k$ on $\Omega$.  By compactness and the continuity
of $X$, there exist
constant $M$ and a non-increasing continuous function $C:\R_+ \ra \R_+$ with $C(0)=0$ so that 
\[ M = \max_\Omega |X|\]
and for $x,y \in \Omega$, 
\[|X(x)-X(y)| \le C(|x-y|)\]
By the compactness of $\Omega$ and the continuity of $X_k$, we have,
for each $k$, constants $M(k)$ and non-increasing continuous functions
$C_k:\R_+ \ra \R_+$ with $C_k(0)=0$ so that 
\[ M(k) = \max_\Omega |X_k|\]

Moreover, since $X_k \ra X$ uniformly on $\Omega$, there exists a
constant $\alpha$ so that 
\[ M(k) \le \alpha M\]

for all $k$.
Next, we construct integral curves for $X$ and $X_k$ emanating from a
point $x_0 \in \Omega$ using Picard's approximation
method.  To do so, let 
\[\phi^0_0(t)=0, \; \; \phi^k_0(t)=0\]
and
\[\phi^0_n(t)= x_0 + \int_0^t X(\phi^0_{n-1}(s))\; ds, \;\; \phi^k_n(t)= x_0 + \int_0^t X_k(\phi^k_{n-1}(s))\; ds \]

\begin{Lem} $\{\phi^0_n\}$ has a subsequence which converges uniformly
  on $\Omega$.
\end{Lem}
\pf  First, since 
\begin{equation*}
\begin{split}
|\phi^0_n(t)| &= \left |x_0+  \int_0^t X(\phi^0_{n-1}(s))\; ds\right| \\
&\le |x_0|+ \int_0^t |X(\phi_{n-1}^0(s))| \; ds \\
&\le |x_0|+ M t
\end{split}
\end{equation*}
we have that the sequence is pointwise bounded.  Second, since
\begin{equation*}
\begin{split}
|\phi^0_n(t_1) - \phi^0_n(t_0)| &= \left |\int_{t_0}^{t_1} X(\phi_{n-1}(s))
 \; ds \right | \\
& \le M|t_1 -t_0|
\end{split}
\end{equation*}
we have that the sequence is equicontinuous (in fact uniformly
Lipschitz).  By the theorem of Arzela-Ascoli, there exists a
subsequence that converges uniformly on $\Omega$.  $\qed$

We note that the same argument applies (with appropriately defined
constants) for $\{\phi^k_n\}$ as well:
\begin{Lem}\label{l2} $\{ \phi^k_n\}_{k,n}$ has a subsequence which converges
  uniformly (in both $k$ and $n$) on $\Omega$.
\end{Lem}
\pf As in the previous lemma, we have
\begin{equation*}
\begin{split}
|\phi^k_n(t)| &\le |x_0| + \int_0^t |X_k(\phi^k_{n-1}(s))| \; ds \\
&\le |x_0| +t M(k) \\
&\le |x_0| + t \alpha M
\end{split}
\end{equation*}
In other words, the sequence is pointwise bounded in both $k$ and
$n$.  Next, 
\begin{equation*}
\begin{split}
|\phi^k_n(t_1)-\phi^k_n(t_0)| & \le
 \int_{t_0}^{t_1}|X_k(\phi^k_{n-1}(s))| \; ds \\
& \le M(k) |t_1-t_0| \\
& \le \alpha M |t_1-t_0|
\end{split}
\end{equation*}
And so the sequence is equicontinuous as well.  Thus, by the theorem of Arzela-Ascoli, we may extract a subsequence
that converges uniformly in both $k$ and $n$ on $\Omega$.  $\qed$

For the purposes of this
discussion, we assume that we have taken the appropriate subsequences
so that $\phi^0_n \ra \phi^0$ and $\phi^k_n \ra \phi^k$ uniformly on
$\Omega$.  This gives us the existence of integral curves for these
vector fields.  Of course, these integral curves may not be unique.    

We next show that, using these integral curves, $X_k(\phi^k(t)) \ra
X(\phi^0(t))$.

\begin{Lem} \[ |\phi^k_n(t)-\phi^0_n(t)| \le M_k t +C^{n-1}(k,t)\]
where 
\[C^m(k,t)= t C(M_kt + t C(M_k t + t C(M_k t + \dots +tC(M_k t))))\]
and the nested applications of the function $C$ occur $m$ times.
\end{Lem}
\pf  We proceed by induction.  First, we note several initial cases:

\[ |\phi^k_0(t)-\phi^0_0(t)| =0 \]
\[ |\phi^k_1(t)-\phi^0_1(t)| \le \int_0^t |X_k(0)-X(0)| \; ds \le M_k t = M_kt + C^0(k,t) \]

\begin{equation}
\begin{split}
|\phi^k_2(t)-\phi^0_2(t)| &\le \int_0^t
 |X_k(\phi^k_1(s))-X(\phi^0_1(s))| \; ds \\
&\le M_k t + \int_0^t |X(\phi^k_1(s))-X(\phi^0_1(s))| \; ds \\
&\le M_k t +\int_0^tC(|\phi^k_1(s)-\phi^0_1(s)|) \; ds \\
&\le M_kt + t C(M_k t)  \;\;\; \text{(by the previous calculuation)}\\
&= M_kt + C^1(k,t)\\
\end{split}
\end{equation}

Now, assume that 
\[|\phi^k_{n-1}(t)-\phi^0_{n-1}(t)| \le M_k t +t C^{n-2}(k,t)\]
Then,
\begin{equation}
\begin{split}
 |\phi^k_n(t)-\phi^0_n(t)| &\le \int_0^t |X_k(\phi^k_{n-1}(s))-X(\phi^0_{n-1}(s))| \; ds \\
&\le M_k t + \int_0^t |X(\phi^k_{n-1}(s))-X(\phi^0_{n-1}(s))| \; ds\\
&\le M_k t +\int_0^tC(|\phi^k_{n-1}(s)-\phi^0_{n-1}(s)|) \; ds \\
&\le M_kt +  \int_0^t C(M_k t + C^{n-1}(k,t)) \; ds  \;\;\; \text{(by the induction hypothesis)}\\
&\le M_k t + tC((M_k t + C^{n-1}(k,t))) = M_k t + C^n(k,t)\\
\end{split}
\end{equation}
This completes the induction and the proof.  
$\qed$

We note that as $M_k$ is a coefficient in every term of each argument
of $C$ in $C^m(k,t)$ and $C(0)=0$, we have that $\lim_{k \ra \infty}
C^m(k,t)=0$ as $\lim_{k \ra \infty}M_k =0$.  Moreover, in light of lemma \ref{l2}, we know that
$\phi^k_n$ tends to some function uniformly as $k \ra \infty$, we see
that the previous lemma implies that $\phi^k_n \ra \phi^0_n$ as $k \ra
\infty$.  

We now prove the claim:
\begin{Lem}\label{l3} \[\lim_{k \ra \infty} |X_k(\phi^k(t))-X(\phi^0(t))| =0\]
\end{Lem}
\pf
\begin{equation*}
\begin{split}
\lim_{k \ra \infty}|X_k(\phi^k(t))-X(\phi^0(t))| &=  \lim_{k \ra
  \infty} \lim_{n \ra \infty} |X_k(\phi^k_n(t))-X(\phi^0_n(t))|\\
 &= \lim_{n \ra
  \infty} \lim_{k \ra \infty} |X_k(\phi^k_n(t))-X(\phi^0_n(t))|\\
&= \lim_{n \ra \infty} |X(\phi^0_n(t)) - X(\phi^0_n(t))| =0
\end{split}
\end{equation*}

In the second equality, we may switch the order of the limits because
$\phi^k_n \ra \phi^k$ uniformly in both $\Omega$ and $k$ as $n \ra 0$
by lemma \ref{l2}.  $\qed$


\end{document}